\documentclass[12pt]{amsart}
\usepackage{amssymb,epsfig,times}

\newtheorem{Thm}{Theorem}
\newtheorem{Cor}[Thm]{Corollary}
\newtheorem{Lemma}[Thm]{Lemma}
\newtheorem{Prop}[Thm]{Proposition}

\theoremstyle{definition}
\newtheorem{Defn}{Definition}

\newtheorem{Remark}{Remark}
\newtheorem{Ex}[Remark]{Example}

\newcommand{\mf}[1]{\mathbb{#1}}
\newcommand{\mc}[1]{\mathcal{#1}}
\newcommand{\mb}[1]{\mathbf{#1}}

\DeclareMathOperator{\Part}{\mathcal{P}}
\DeclareMathOperator{\NC}{\mathit{NC}}
\DeclareMathOperator{\Sym}{Sym}
\DeclareMathOperator{\E}{\mathrm{E}}
\DeclareMathOperator{\Int}{\mathit{Int}}
\DeclareMathOperator{\Outer}{\mathrm{Outer}}
\DeclareMathOperator{\Sing}{\mathrm{Sing}}
\DeclareMathOperator{\cyc}{\mathrm{cyc}}

\newcommand{\norm}[1]{\left\Vert#1\right\Vert}
\newcommand{\abs}[1]{\left\vert#1\right\vert}
\newcommand{\chf}[1]{\mathbf{1}_{#1}}
\newcommand{\set}[1]{\left\{#1\right\}}
\newcommand{\ip}[2]{\left \langle #1, #2 \right \rangle}
\newcommand{\state}[1]{\varphi \left[ #1 \right]}
\newcommand{\rc}[1]{\mathrm{rc} \left( #1 \right)}
\newcommand{\Exp}[1]{\left\langle #1 \right\rangle}
\newcommand{\Cexp}[2]{\E_{#1} \left[ #2 \right]}
\newcommand{\A}[1]{A \left( #1 \right)}
\newcommand{\Aa}[2]{A_{#1} \left( #2 \right)}
\newcommand{\KS}[1]{\mathrm{W} \left( #1 \right)}

\renewcommand{\phi}{\varphi}
\newcommand{\br}{\medskip\noindent}

\allowdisplaybreaks[1]

\title{Appell polynomials and their relatives}
\author[M.~Anshelevich]{Michael Anshelevich}
\thanks{This work was supported in part by an NSF postdoctoral fellowship}
\address{Department of Mathematics, University of California, Riverside, CA 92521-0135}
\email{manshel@math.ucr.edu}
\date{\today}

\begin{document}


\begin{abstract}
This paper summarizes some known results about Appell polynomials and investigates their various analogs. The primary of these are the free Appell polynomials. In the multivariate case, they can be considered as natural analogs of the Appell polynomials among polynomials in non-commuting variables. They also fit well into the framework of free probability. For the free Appell polynomials, a number of combinatorial and ``diagram'' formulas are proven, such as the formulas for their linearization coefficients. An explicit formula for their generating function is obtained. These polynomials are also martingales for free L\'{e}vy processes. For more general free Sheffer families, a necessary condition for pseudo-orthogonality is given. Another family investigated are the Kailath-Segall polynomials. These are multivariate polynomials, which share with the Appell polynomials nice combinatorial properties, but are always orthogonal. Their origins lie in the Fock space representations, or in the theory of multiple stochastic integrals. Diagram formulas are proven for these polynomials as well, even in the $q$-deformed case.
\end{abstract}

\maketitle

\section{Introduction}
\noindent
Let $\mu$ be a probability measure on the real line all of whose moments
\[
m_n(\mu) = \int_{\mf{R}} x^n \,d\mu(x)
\]
are finite. Then there are at least three natural families of polynomials associated to such a measure. The most familiar ones are the orthogonal polynomials $\set{P_n}_{n=0}^\infty$. This is a polynomial family (that is, $P_n$ has degree $n$) such that
\[
\ip{P_n}{P_k}_{\mu} = \int_{\mf{R}} P_n(x) P_k(x) \,d\mu(x) = 0
\]
for $n \neq k$. Two standard normalizations are to require the polynomials to be monic or to be orthonormal. By the spectral theorem for orthogonal polynomials, such (monic) polynomials satisfy a three-term recursion relation
\[
P_{n+1}(x) = x P_n(x) - \alpha_n P_n(x) - \beta_n P_{n-1}(x),
\]
where $\alpha_n \in \mf{R}$, $\beta_n \in \mf{R}_+$ are the Jacobi parameters, and $P_{-1} = 0, P_0(x) = 1$. The Jacobi parameters and the moments of the measure can be expressed in terms of each other, and their properties related to the properties of the measure and the orthogonal polynomials, for example using the Viennot-Flajolet theory \cite{Flajolet,Viennot-Notes,Viennot-Short}. A typical question in this direction is to find explicitly the \emph{linearization coefficients}
\[
\int_{\mf{R}} P_{n_1}(x) P_{n_2}(x) \ldots P_{n_k}(x) \,d\mu(x).
\]
Already the proofs of the positivity of these coefficients are quite subtle \cite{deM-Stanton-Expansions} and they are known explicitly only in very rare cases \cite{KimZeng}.

\br
Another natural and very classical \cite{Appell} polynomial family associated to $\mu$ is its family of Appell polynomials, which have the exponential generating function
\begin{equation}
\label{Appell-Generating}
\sum_{n=0}^\infty \frac{1}{n!} A^{(n)}(x) z^n = \exp(x z - \log M_\mu(z)),
\end{equation}
where
\[
M_\mu(z) = \sum_{n=0}^\infty \frac{1}{n!} m_n(\mu) z^n
\]
is the exponential moment generating function of $\mu$. It is easy to see that an equivalent definition is via a recursion relation
\begin{equation}
\label{Appell-Recursion}
A^{(n+1)}(x) = x A^{(n)}(x) - \sum_{k=0}^n \binom{n}{k} r_{n+1-k}(\mu) A^{(k)}(x)
\end{equation}
and $A^{(0)}(x) = 1$. Here, $r_k(\mu)$ are the cumulants (semi-invariants) of the measure. Unlike the Jacobi parameters, the cumulants are obtained from the moments of the measure via a simple relation
\[
\sum_{k=1}^\infty \frac{1}{k!} r_k(\mu) z^k = \log M_{\mu}(z)
\]
(both sides considered as formal power series). Many classical (non-orthogonal) polynomial families are Appell. They arise in finite operator calculus \cite{RotaFiniteCalculusBook} and the study of hypergroups, in numerical analysis (Bernoulli polynomials are Appell), but also in probability theory, in the study of stochastic processes \cite{Lai}, non-central limit theorems \cite{Avram-Taqqu,Giraitis-Surgailis}, and natural exponential families \cite{Pommeret-Appell}. From the combinatorial point of view, they have nice linearization and multinomial formulas.

\br
The third family of polynomials has not apparently been explicitly defined before, although it appears implicitly in the paper \cite{Kailath-Segall}. For this reason, we will call them the Kailath-Segall polynomials. These are polynomials in (infinitely many) variables $\set{x_k}_{k=1}^\infty$. They are indexed by all finite sequences of natural numbers $\vec{u} = (u(1), u(2), \ldots, u(n))$, $n \geq 0$, and defined by the recursion
\begin{multline}
\label{Multi-KS}
W_{j, u(1), u(2), \ldots, u(n)} = x_j W_{u(1), u(2), \ldots, u(n)} \\
- \sum_{i=1}^n r_{j + u(i)}(\mu) W_{u(1), \ldots, \widehat{u(i)}, \ldots, u(n)} - \sum_{i=1}^n W_{j + u(i), u(1), \ldots, \widehat{u(i)}, \ldots, u(n)} - r_j(\mu) W_{u(1), u(2), \ldots, u(n)}
\end{multline}
with initial conditions $W_{\emptyset} = 1$, $W_i = x_i - r_i$. As usual, $\widehat{u(i)}$ means ``omit the $i$'th term''. Note that $W_{\vec{u}}$ contains a single monomial of the highest degree $\abs{\vec{u}} = \sum_{i=1}^n u(i)$, namely $x_{\vec{u}} = x_{u(1)} x_{u(2)} \ldots x_{u(n)}$, and that it is a polynomial in the variables
\[
\set{x_i: i = \sum_{j \in S} u(j) \text{ for some } S \subset \set{1,2, \ldots, n}}.
\]
The origin of these polynomials is again in probability theory, where they appear as certain multiple stochastic integrals. Moreover, they share a number of properties with both the orthogonal and the Appell polynomials. Their recursion relation is determined by the cumulants, and they have nice linearization properties, just like the Appell polynomials. In fact, we will show that $A^{(n)}(x_1)$ is a linear combination of the Kailath-Segall polynomials. On the other hand, the polynomials $W^{(n)}= W_{1, 1, \ldots, 1}$ also have the following orthogonality property. Define a measure $\mu^{(n)}$ on $\mf{R}^n$ by specifying its multivariate cumulants:
\[
r_{\vec{u}}(\mu^{(n)}) = r_{\abs{\vec{u}}} (\mu).
\]
See Remark~\ref{Remark:measure} for the fashion in which the measure $\mu^{(n)}$ is determined by its cumulants. If $\mu$ is infinitely divisible, $\mu^{(n)}$ is a positive measure. Then
\[
\int_{\mf{R}^n} W^{(n)}(x_1, x_2, \ldots, x_n) W^{(k)}(x_1, x_2, \ldots, x_k) \,d\mu^{(n)}(x_1, x_2, \ldots, x_n) = 0
\]
if $k < n$.

\br
A generalization of the Appell polynomials are the Sheffer polynomials. Let $U$ be a function with a formal power series expansion such that $U(z) = z + \text{ higher-order terms}$. Then the Sheffer polynomials are defined via their exponential generating function
\[
\sum_{n=0}^\infty \frac{1}{n!} P_n(x) z^n = \exp \bigl(x U(z) - \log M_{\mu}(U(z)) \bigr).
\]
Equivalently,
\[
U^{-1}(\partial_x)(P_n) = n P_{n-1} \text{ and } \ip{P_n}{1}_\mu = 0 \text{ for } n>0.
\]
They share with the Appell polynomials the multinomial expansion properties and the relation to stochastic processes (see Section~\ref{Subsec:Martingales}). Among the Appell polynomials, only the Hermite ones are orthogonal. Meixner's classic characterization \cite{Meixner} describes all the orthogonal Sheffer polynomials. There are also multivariate versions of this statement \cite{Pommeret-Sheffer}.

\br
We start the paper by describing some properties of and relations between the three aforementioned families of polynomials. In the Appell and Kailath-Segall cases, natural starting points are in fact certain multi-linear functionals on more general algebras, which can then be specialized to polynomials. These families have similar but different properties. For example, the Appell linearization coefficients are sums over non-homogeneous partitions, while the Kailath-Segall linearization coefficients are sums over inhomogeneous partitions. A number of the results in Section~\ref{Section:Classical} are known and so are stated without proof.

\br
As mentioned above and described in more detail below, both the Appell and the Kailath-Segall polynomials arise in probability theory and are related to the notion of \emph{independence}. Let $\phi$ be a real linear functional on the algebra $\mf{R}[x, y]$ of polynomials in two variables. Then $x, y$ are independent with respect to $\phi$ if for any $P, Q$,
\[
\state{P(x) Q(y)} = \state{P(x)} \state{Q(y)}.
\]
In the early 1980's, Dan Voiculescu introduced a parallel but really very different notion of \emph{free independence} \cite{Voi85}. Let $\phi$ now be a real linear functional on the algebra $\mf{R} \langle x, y \rangle$ of polynomials in two \emph{non-commuting} variables. Then $x, y$ are freely independent with respect to $\phi$ if whenever
\[
\state{P_1(x)} = \state{Q_1(y)} = \ldots = \state{P_n(x)} = \state{Q_n(y)} = 0
\]
and $Q_0, P_{n+1}$ each are either centered or scalar, then
\[
\state{Q_0(y) P_1(x) Q_1(y) \ldots P_n(x) Q_n(y) P_{n+1}(x)} = 0.
\]
A whole theory of \emph{free probability} \cite{VoiSF}, based on this notion, is by now quite well developed. It turns out that there are ``free analogs of'' the Appell and Kailath-Segall polynomials, which, very roughly, are obtained by replacing commuting variables with non-commuting ones, exponential generating functions with the usual ones, and exponentials with resolvents. Such a replacement, however, is quite non-trivial. The analysis of the preceding section is repeated in Section~\ref{Section:Free} for the free case, except that in this case most of the results are new. For the free Appell polynomials, we find an explicit form of the generating function and various recursion and ``diagram'' formulas. As expected, these formulas are in most cases based on the lattice of non-crossing, rather than all, set partitions (but not always, compare for example equation~\eqref{Appell-X} with Proposition~\ref{Prop:Free-Appell}(c)). For the free Sheffer polynomials, we find a necessary condition for them to be pseudo-orthogonal. Consequences of these results for free probability will be developed elsewhere.

\br
Comparison of formulas from the preceding two sections shows that many of them appear as particular cases of $q$-interpolated forms, with the usual case corresponding to $q=1$ and the free case corresponding to $q=0$. On the other hand, many other formulas do not appear to admit of such an interpolation. In Section~\ref{Section:Q}, we define the $q$-Appell and $q$-Kailath-Segall polynomials, and show that some of their properties carry over to the whole interpolated family. Some of the other properties, at least at present, do not. Therefore this section is necessarily more tentative than the preceding ones. We only consider single-variable $q$-Appell polynomials, and find an explicit form of the generating function for them, as well as the $q$-analogs of the Meixner families. For the $q$-Kailath-Segall polynomials, we find various recursion and ``diagram'' formulas. Finally, in the appendix we show that $q$-Appell polynomials are not linear combinations of the $q$-Kailath-Segall polynomials, unlike in the classical and the free case. As a consequence, they cannot be martingale polynomials for the $q$-L\'{e}vy processes.

\br
\noindent\textbf{Acknowledgements:} Thanks to Mourad Ismail, Marius Junge, and Michael Skeide for useful and enjoyable conversations. Thanks also to the referees for a number of helpful suggestions and criticisms, especially the correction in Definition~\ref{Defn:Multi-free-Appell}.

\section{Classical polynomial families}
\label{Section:Classical}

\subsection{Notation}
We will use multi-index notation $\vec{u} = (u(1), u(2), \ldots, u(k))$. Denote by
\[
\Delta_{k,n} = \set{\vec{u} \in \mf{N}^k: \abs{\vec{u}} = n}
\]
a basic simplex (where $0 \not \in \mf{N}$). For two multi-indices $\vec{u}, \vec{v}$, $(\vec{u}, \vec{v})$ will denote their concatenation. For $\vec{u} \in \mf{N}^n$, $B \subset \set{1, 2, \ldots, n}$, $B = \set{v(1) < v(2) < \ldots < v(k)}$, denote
\[
(\vec{u}:B) = (u(v(1)), u(v(2)), \ldots, u(v(k))).
\]
For a subset $A \subset B$, $A^c$ will denote the complement of $A$, where $B$ is understood.

\subsubsection{Partitions.}
A set partition of a set $S$ is a collection of disjoint non-empty subsets of $S$ whose union is $S$. If $S$ is an ordered set and $\pi = \set{B_1, B_2, \ldots, B_k}$ is such a partition, we order the classes of $\pi$ according to the order of their smallest elements, $\min(B_1) < \min(B_2) < \ldots  < \min(B_k)$. We will consider the following three lattices of partitions. By $\Part(n)$ we'll denote the lattice of all partitions of the set $\set{1, 2, \ldots, n}$. For $\pi \in \Part(n)$, $\Sing(\pi)$ is the collection of single-element (singleton) classes of $\pi$.

\br
By $\NC(n)$ we'll denote the lattice of \emph{non-crossing partitions} \cite{Kre72}. These are the partitions with the property that
\[
i < j < k, \; i \stackrel{\pi}{\sim} k , \; j \stackrel{\pi}{\sim} l,  \; i \stackrel{\pi}{\not \sim} j \;
\Rightarrow \; i < l < k.
\]
For a non-crossing partition $\pi$, a class $B$ is called \emph{outer} if for any other class $C \in \pi$, if $i, j \in C$, $k \in B$, and $i < k$, then $j < k$. Otherwise a class is called inner. The outer classes of $\pi$ will be denoted by $\Outer(\pi)$.

\br
The third lattice of partitions, which we use mostly for notational convenience, is that of \emph{interval partitions}, all of whose classes are intervals of consecutive integers. This lattice $\Int(n)$ is naturally isomorphic to the power set of $\set{1, 2, \ldots, n-1}$. For $\sum_{i=1}^k s(i) = n$, we denote by $\pi_{s(1), s(2), \ldots, s(k)} \in \Int(n)$ the interval partition
\begin{equation*}
\set{\Bigl(1, \ldots, s(1) \Bigr), \Bigl(s(1)+1, \ldots, s(1)+s(2) \Bigr), \ldots, \Bigl( \sum_{i=1}^{k-1} s(i) + 1, \ldots, n \Bigr)}.
\end{equation*}

\br
There is a partial order $\leq$ on $\Part(n)$ which restricts to the other two lattices. We denote the smallest element in $\Part(n)$ by $\hat{0} = \set{(1), (2), \ldots, (n)}$ and the largest one by $\hat{1} = \set{(1, 2, \ldots, n)}$. We denote the meet and the join in the lattices by $\wedge$ and $\vee$, respectively. In particular,
\[
i \stackrel{\pi \wedge \sigma}{\sim} j  \Leftrightarrow i \stackrel{\pi}{\sim} j \text{ and } i \stackrel{\sigma}{\sim} j.
\]

\begin{figure}[ht]
  \psfig{figure=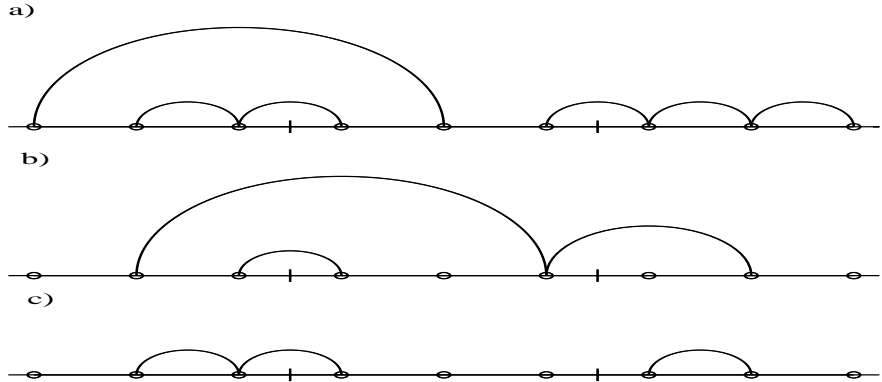,height=2in,width=0.7\textwidth}
  \caption{With respect to the partition $\set{(1,2,3), (4,5,6), (7,8,9)}$, partitions which are (a) non-homogeneous, connected, (b) inhomogeneous, connected, (c) not connected.}
  \label{Figure:Homogeneous}
\end{figure}

\begin{Defn}
Let $\sigma \in \Part(N)$ be a partition. For a partition $\pi \in \Part(N)$, we say that
\begin{enumerate}
\item
A class $B \in \pi$ is \emph{homogeneous} with respect to $\sigma$ if $B \subset C$ for some $C \in \sigma$,
\item
$\pi$ is \emph{non-homogeneous} with respect to $\sigma$ if $\pi$ has no homogeneous classes with respect to $\sigma$,
\item
$\pi$ is \emph{inhomogeneous} with respect to $\sigma$ if $\pi \wedge \sigma = \hat{0}$, that is, $i \stackrel{\pi}{\sim} j \Rightarrow i \stackrel{\sigma}{\not \sim} j$,
\item
$\pi$ is \emph{connected} with respect to $\sigma$ if $\pi \vee \sigma = \hat{1}$.
\end{enumerate}
See Figure~\ref{Figure:Homogeneous}.
\end{Defn}

\subsubsection{Extended partitions and restricted crossings}
\begin{figure}[ht]
  \psfig{figure=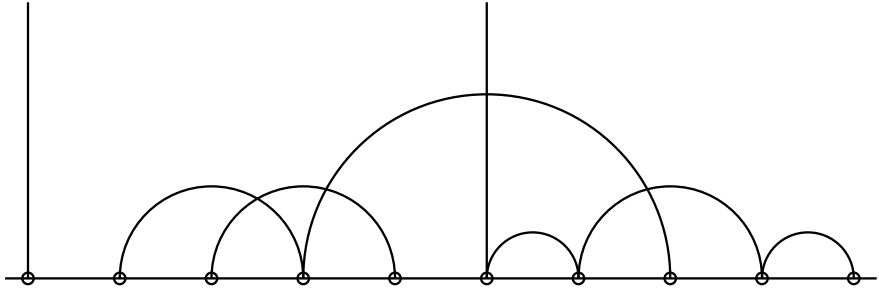,height=1.5in,width=0.7\textwidth}
  \caption{An extended partition of $10$ elements with $2$ left-open classes and $4$ restricted crossings.}
  \label{Figure:Extended}
\end{figure}

For $S \subset \pi$, call the pair $(S, \pi)$ an extended partition; $S$ is to be thought of as the collection of classes of $\pi$ ``open on the left''. See Figure~\ref{Figure:Extended} for an example.

\noindent
For $1 \leq k \leq m \leq n$, define the restriction
\[
(S', \pi') = (S, \pi) \upharpoonright \set{k, \ldots, m}
\]
as follows:
\begin{align*}
B' \in \pi' &\text{ if } B' = B \cap \set{k, \ldots, m}, B \in \pi, \\
B' \in S' &\text{ if } B \in S \text{ or } B \cap \set{1, \ldots, k-1} \neq \emptyset.
\end{align*}
Define the number of right restricted crossings of $(S, \pi)$ at the point $k$ as follows:
\[
\rc{k, S, \pi} =
\begin{cases}
0, & \text{if } k \in B, k = \max(B) \text{ or } k+1 \in B, \\
\abs{S'}, & \text{if } k \in B, j = \min \set{i \in B, i > k}, \\
& \qquad (S', \pi') = (S, \pi) \upharpoonright \set{k+1, \ldots, j-1}.
\end{cases}
\]
Let $\rc{S, \pi} = \sum_{k=1}^n \rc{k, S, \pi}$. Note that also
\[
\rc{S, \pi}
= \rc{\pi} + \sum_{B \in S} \abs{C \in \pi: \min(C) < \min(B) < \max(C)},
\]
where $\rc{\pi} = \rc{\emptyset, \pi}$ (see \cite{BiaCross}).

\subsubsection{Cumulants}
A measure $\mu$ on $\mf{R}$ all of whose moments are finite induces a positive semi-definite unital linear functional $\phi$ on $\mf{R}[x]$ by $\state{x^n} = m_n(\mu)$. Positivity will not play a part in most of the results below. Thus, our starting object is a linear functional $\phi$, which does not necessarily correspond to a positive measure, although we still assume that it is unital, $\state{1} = 1$. Throughout the paper, the functional will be fixed, and so will frequently be omitted from notation.

\br
In the multi-dimensional situation, let $\phi$ be  a general unital real linear functional on some real algebra $\mc{A}$. For $X_1, X_2, \ldots, X_n \in \mc{A}$, denote by
\[
M[X_1, X_2, \ldots, X_n]
= \state{X_1 X_2 \ldots X_n}
\]
the \emph{joint moment} of $(X_1, X_2, \ldots, X_n)$ with respect to $\phi$. Also, for $\pi \in \Part(n)$, denote
\[
M_{\pi}[X_1, X_2, \ldots, X_n]
= \prod_{B \in \pi} M[X_i: i \in B]
\]
the partitioned moment of $(X_1, X_2, \ldots, X_n)$. Here, in the case when $\mc{A}$ is non-commutative, the factors in $\prod_{i \in B} X_i$ and the terms in $\set{X_i: i \in B}$ are taken in order. For each $\pi$, $M_{\pi}$ is an $n$-linear map. Given $M$, define the corresponding \emph{joint cumulant} recursively by
\begin{equation}
\label{Moment-cumulant}
R[X_1, X_2, \ldots, X_n] = M[X_1, X_2, \ldots, X_n] - \sum_{\substack{\pi \in \Part(n) \\ \pi \neq \hat{1}}} R_{\pi}[X_1, X_2, \ldots, X_n],
\end{equation}
where the partitioned cumulants $R_\pi$ are defined as above, and involve only cumulants of lower order. These are also multi-linear maps.

\br
The multi-index notation will be used extensively but consistently. Throughout the paper
\begin{gather}
\notag
X_{\vec{u}} = X_{u(1)} X_{u(2)} \ldots X_{u(k)}, \\
\label{Moment-consistent}
M[X_{\vec{u}}] = M[X_{u(1)}, X_{u(2)}, \ldots, X_{u(k)}] = M[X_{u(1)} X_{u(2)} \ldots X_{u(k)}], \\
\notag
R[X_{\vec{u}}] = R[X_{u(1)}, X_{u(2)}, \ldots, X_{u(k)}], \\
\notag
R[(X_{\vec{u}})] = R[X_{u(1)} X_{u(2)} \ldots X_{u(k)}]
\end{gather}
On the other hand, $\mb{x} = (x_1, x_2, \ldots, x_n)$ will denote the total collection of variables involved, and same for $z, w, X$, etc. Finally, for $B \subset \set{1, \ldots, n}$,
\[
x_B = \prod_{i \in B} x_i,
\]
the product taken in order of increasing indices.

\br
For $S = \set{B_1, \ldots, B_k}$ a collection of disjoint subsets of $\set{1, \ldots, n}$ (for example, a partition), $S$ is always ordered according to the order of the smallest elements of the subsets.

\begin{Remark}
\label{Remark:measure}
If $\mu$ is a measure on $\mf{R}^n$ all of whose moments
\[
m_{\vec{u}}(\mu) = \int_{\mf{R}^n} x_{\vec{u}} \,d\mu(\mb{x})
\]
are finite, its cumulants $r_{\vec{u}}(\mu)$ are defined in terms of the moments by equation~\eqref{Moment-cumulant}. Note that while the functional $x_{\vec{u}} \mapsto m_{\vec{u}}(\mu)$ is determined by its cumulants, again via equation~\eqref{Moment-cumulant}, the measure $\mu$ itself may not be. The determinacy of the moment problem in the multivariate context is a difficult question, see \cite[Chapter 3]{Dunkl-Xu} and their references.
\end{Remark}

\subsection{Appell polynomials}
Appell polynomials are defined by equations~\eqref{Appell-Generating} or \eqref{Appell-Recursion}. They satisfy the following properties, easily obtained using the generating function:
\begin{gather}
\label{Appell-centered}
\partial_x A^{(n)}(x) = n A^{(n-1)}(x), \quad \state{A^{(n)}} = \delta_{n 0}, \\
\label{x-Appell}
x^n = \sum_{k=0}^n \binom{n}{k} m_{n-k} A^{(k)}(x), \\
\label{Appell-x}
A^{(n)}(x) = \sum_{k=0}^n (m^{-1})_{n-k} x^k,
\end{gather}
where
\[
(m^{-1})_n = \sum_{k=1}^n \sum_{\vec{u} \in \Delta_{k, n}} (-1)^k \frac{1}{k!} \binom{n}{u(1), \ldots, u(k)} r_{u(1)} \ldots r_{u(k)},
\]
and the notation is suggested by the relation
\[
\sum_{k=0}^n \binom{n}{k} m_k (m^{-1})_{n-k} = \delta_{n 0}.
\]
See the original paper \cite{Appell}, where \eqref{Appell-centered} was taken as the definition, or \cite{Giraitis-Surgailis}.

\begin{Remark}
The following result of Appell~\cite{Appell}, developed in great detail by Rota et al.~\cite{RotaFiniteCalculusBook}, appears to have no free analog. Let $\set{A^{(n)}}$, $\set{B^{(n)}}$ be two families of Appell polynomials, with exponential generating functions $F, G$, respectively. Define a new family of polynomials $\set{(A B)^{(n)}}$ as follows: expand $A^{(n)}$ in the powers of $x$, and for each $x^k$ substitute $B^{(k)}$. Then $(A B)^{(n)} = (B A)^{(n)}$, and these polynomials are again an Appell family, with exponential generating function $F G$.
\end{Remark}

\subsection{Multivariate Appell polynomials \cite{Avram-Taqqu,Giraitis-Surgailis}}
Let $\mc{A}$ be a commutative real algebra with a unital real linear functional $\phi$. We will usually call elements of such an algebra \emph{random variables}, since any collection of real-valued random variables on some probability space, such that all of their joint moments are finite, generates a commutative real algebra with the expectation functional on it. For $n \geq 1$, define an $n$-linear map $A: \mc{A}^n \rightarrow \mc{A}$ as follows.  For $X_1, X_2, \ldots, X_n \in \mc{A}$, $\A{X_1, X_2, \ldots, X_n}$ is in fact a symmetric polynomial in $\set{X_i}_{i=1}^n$, which we denote by 
\[
\Aa{X_1, X_2, \ldots, X_n}{x_1, x_2, \ldots, x_n}.
\]
It is determined recursively by
\begin{equation}
\label{Classical-Delta-An}
\partial_{x_i} \Aa{X_1, X_2, \ldots, X_n}{x_1, x_2, \ldots, x_n} = \Aa{X_1, X_2, \ldots, \hat{X}_i, \ldots, X_n}{x_1, \ldots, \hat{x}_i, \ldots, x_n}
\end{equation}
and
\begin{equation}
\label{Classical-multi-centered}
\state{\A{X_1, X_2, \ldots, X_n}} = 0
\end{equation}
for $n > 0$, where $A_{\emptyset} = 1$ is understood. The advantage of this notation is that we can consider only the algebra generated by $\set{X_i}_{i=1}^n$ (and $1$), and the restriction of $\phi$ to this algebra can be thought of as the joint distribution of these random variables. For $\vec{u} \in \set{1, 2, \ldots, n}^k$, define the polynomials $\Aa{\vec{u}}{\mb{x}}$ by
\[
\Aa{\vec{u}}{\mb{x}}
= \Aa{X_{u(1)}, X_{u(2)}, \ldots, X_{u(k)}}{x_{u(1)}, x_{u(2)}, \ldots, x_{u(k)}}.
\]
Thus
\[
\Aa{\vec{u}}{\mb{X}}
= \A{X_{u(1)}, X_{u(2)}, \ldots, X_{u(k)}}.
\]
Note that this notation differs from the usual one: ${A}_{\vec{u}}$ depends in fact only on the number of occurrences of each index $u(i)$, and not on their order; in the usual notation one writes down the number of such occurrences. Our notation is better suited to the non-commutative case.

\br
The following are some properties of the multivariate Appell polynomials. Denote
\[
\mb{x} \cdot \mb{z} = \sum_{i=1}^n x_i z_i,
\]
and denote by
\[
\mb{R}(z) = \sum_{k=1}^\infty \sum_{\vec{v} \in \set{1, \ldots, n}^k} \frac{1}{k!} R[X_{\vec{v}}] z_{\vec{v}}
\]
the exponential cumulant generating function. Note that the $z_i$'s commute and $R$ is symmetric in its arguments, so each term
\[
z_1^{i(1)} z_2^{i(2)} \ldots z_n^{i(n)}, \quad \sum_{j=1}^n i(j) = k
\]
appears
\[
\binom{k}{i(1), i(2), \ldots, i(n)} = \frac{k!}{i(1)! i(2)! \ldots i(n)!}
\]
times. So this notation coincides with the usual one.

\br
Single- and multivariate Appell polynomials are related via
\[
A^{(n)}(x) = \Aa{1, 1, \ldots, 1}{x, x, \ldots, x},
\]
where the single-variable polynomial corresponds to the moment sequence $m_n = \state{X_1^n}$. Also,
\begin{gather}
\label{Classical-multi-Generating}
1 + \sum_{k=1}^\infty \frac{1}{k!} \sum_{\vec{v} \in \set{1, \ldots, n}^k} \Aa{\vec{v}}{\mb{x}} z_{\vec{v}} = \exp \bigl( \mb{x} \cdot \mb{z} - \mb{R}(\mb{z}) \bigr), \\
\label{Classical-recursion}
\A{X_j, X_{\vec{u}}} = X_j \A{X_{\vec{u}}} - \sum_{V \subset \set{1, \ldots, n}} R[X_j, X_{(\vec{u}:V)}] \A{X_{(\vec{u}:V^c)}}, \\
\notag
X_{\vec{u}} = \sum_{\pi \in \Part(n)} \sum_{B \in \pi} \A{X_{(\vec{u}:B)}} \prod_{\substack{C \in \pi, \\ C \neq B}} R[X_{(\vec{u}:C)}]
= \sum_{V \subset \set{1, \ldots n}} \A{X_{(\vec{u}:V)}} \sum_{\pi \in \Part(V^c)} M_\pi[X_{(\vec{u}:V^c)}], \\
\label{Appell-X}
\A{X_{\vec{u}}} = \sum_{\pi \in \Part(n)} \sum_{B \in \pi} X_{(\vec{u}:B)} (-1)^{\abs{\pi}-1} \prod_{\substack{C \in \pi, \\ C \neq B}} R[X_{(\vec{u}:C)}].
\end{gather}
Here $V^c$ is the complement of $V$ in $\set{1, \ldots, n}$. Equation~\eqref{Classical-multi-Generating} looks slightly unusual because of our different notation. The single-variable polynomials also satisfy a binomial formula: if $X, Y$ are independent,
\begin{equation}
\label{binom}
A^{(n)}(X + Y) = \sum_{k=0}^n \binom{n}{k} A^{(k)}(X) A^{(n-k)}(Y).
\end{equation}
This is, of course, the short-hand notation for
\[
A^{(n)}_{X+Y}(x+y) = \sum_{k=0}^n \binom{n}{k} A^{(k)}_X(x) A^{(n-k)}_Y(y).
\]

\begin{Remark}
For fixed $\set{X_1, X_2, \ldots, X_n} \subset \mc{A}$, their joint distribution is the functional $\phi_{\mb{X}}$ on $\mf{R} \langle x_1, x_2, \ldots, x_n \rangle$ determined by
\[
\phi_{\mb{X}}[x_{\vec{u}}] = \state{X_{\vec{u}}}.
\]
Any such functional $\psi$ on $\mf{R} \langle x_1, x_2, \ldots, x_n \rangle$ can be included in a one-parameter family $\set{\psi_t}_{t \in [0, \infty)}$ in the following fashion. Let $R$ denote the joint cumulants of $\set{x_1, x_2, \ldots, x_n}$ with respect to $\psi$. Then define
\[
\psi_t[x_{\vec{u}}]
= \sum_{\pi \in \Part(n)} t^{\abs{\pi}} \prod_{B \in \pi} R[x_{(\vec{u}:B)}].
\]
That is, the joint moments of $\set{x_1, x_2, \ldots, x_n}$ under $\psi_t$ correspond to the joint cumulants $t R$ via relation~\eqref{Moment-cumulant}. If all the linear functionals $\psi_t$ are positive, the corresponding measures form a semigroup with respect to the convolution operation. They are also marginal distributions of the corresponding L\'{e}vy process (see Section~\ref{Subsec:Martingales}). All $\set{\psi_t}$ are positive if and only if $\psi$ itself is infinitely divisible. Without the positivity requirement, any functional is algebraically infinitely divisible. So any family of Appell polynomials naturally comes included in a one-parameter family $A^{(n)}(\cdot, t)$.

\br
A similar construction, based on relation~\eqref{Free-moment-cumulant} can be done in the free case (see Section~\ref{Section:Free}). In this case, if all the linear functionals $\psi_t$ are positive, the corresponding measures form a semigroup with respect to the additive free convolution operation, and they are marginal distributions of the corresponding free L\'{e}vy process. This is the case if $\psi$ itself is freely infinitely divisible. However, in the free one-dimensional case, more is true: for any positive $\psi$, $\psi_t$ is positive for $t \geq 1$ \cite{Nica-Speicher-Multiplication}.
\end{Remark}

\subsection{Martingales}
\label{Subsec:Martingales}
Let $\set{X(t)}_{t \in [0, \infty)}$ be a L\'{e}vy process, that is, a stochastic process with stationary independent increments. Assume that all of the joint moments of $\set{X(t)}$  (with respect to the expectation functional) are finite. Let $\E_t$ be the conditional expectation onto the (von Neumann) algebra generated by $\set{X(s)}_{s \in [0, t)}$, which extends to an orthogonal projection on the space of all square-integrable random variables. Then for each $n$, $A^{(n)}(X(t))$ is a martingale, that is,
\[
\Cexp{s}{A^{(n)}(X(t))} = A^{(n)}(X(s)).
\]
The following are two elementary proofs of this fact. If $M(z)$ is the moment generating function for $X(1)$, then the moment generating function for $X(t)$ is $M(z)^t$. So using the generating function~\eqref{Appell-Generating} of the Appell polynomials,
\[
\begin{split}
&\Cexp{s}{\exp \bigl(X(t) z - t \log M(z) \bigr)} \\
&\qquad = \Cexp{s}{\exp \Bigl((X(t) - X(s)) z - (t-s) \log M(z) \Bigr) \exp \Bigl(X(s) z - s \log M(z) \Bigr)} \\
&\qquad = \exp(X(s) z - s \log M(z))
\end{split}
\]
since $\E \left[\exp \bigl((X(t) - X(s)) z \bigr) \right] = M(z)^{t-s}$. On the other hand, using the binomial property~\eqref{binom},
\[
\begin{split}
\Cexp{s}{A^{(n)}(X(t))}
&= \Cexp{s}{A^{(n)} \Bigl(X(s) + \bigl(X(t) - X(s)\bigr) \Bigr)} \\
&= \Cexp{s}{\sum_{k=0}^n \binom{n}{k} A^{(k)} \bigl(X(s) \bigr) \cdot A^{(n-k)} \bigl(X(t) - X(s) \bigr)}
= A^{(n)}(X(s)),
\end{split}
\]
since $\E \left[A^{(n-k)}(X(t) - X(s)) \right] = 0$ for $n > k$. Independence of increments and properties of the conditional expectation are used in both proofs.

\begin{Defn}
A \emph{monic polynomial family} in $n$ variables is a family of polynomials indexed by
\[
\set{\vec{u} \in \bigcup_{k=0}^\infty \set{1, 2, \ldots, n}^k}
\]
such that $P_{\emptyset}(\mb{x}) = 1$ and each
\[
P_{\vec{u}}(\mb{x}) = x_{\vec{u}} + \text{ lower order terms},
\]
where the grading is by total degree in $x_1, x_2, \ldots, x_n$. All the families of polynomials considered in this paper form monic polynomial families.

\br
Let $\mf{R}\langle \mb{x} \rangle = \mf{R}\langle x_1, x_2, \ldots, x_n \rangle$ be all the real polynomials in $n$ non-commuting variables. For a multi-index $\vec{u} \in \mf{N}^k$, denote
\[
(\vec{u})^{op} = (u(k), \ldots, u(2), u(1)).
\]
Define an involution on $\mf{R}\langle \mb{x} \rangle$ via an $\mf{R}$-linear extension of
\[
(x_{\vec{u}})^\ast = x_{(\vec{u})^{op}}.
\]
Similarly, define an involution on $\mf{C}\langle \mb{x} \rangle$ via a $\mf{C}$-anti-linear extension of the same relation.

\br
A monic polynomial family $\set{P_{\vec{u}}}$ is pseudo-orthogonal with respect to a functional $\phi$ if
\[
\state{P_{\vec{u}}^\ast P_{\vec{v}}} = 0
\]
whenever $\abs{\vec{u}} \neq \abs{\vec{v}}$. The family is orthogonal if this is the case whenever $\vec{u} \neq \vec{v}$.
\end{Defn}

\subsection{Fock spaces}
\label{Section:Fock}
Let $\set{P_{\vec{u}}(\mb{X})}$ be a monic polynomial family in $n$ variables. Define a functional $\phi$ on $\mf{C} \langle X_1, X_2, \ldots, X_n \rangle$ by $\state{1} = 1$, $\state{P_{\vec{u}}} = 0$ for $\abs{\vec{u}} \geq 1$, and extend $\mf{C}$-linearly. On $\mf{C} \langle x_1, x_2, \ldots, x_n \rangle$, define an inner product via
\[
\ip{x_{\vec{u}}}{x_{\vec{v}}}
= \state{P_{\vec{u}}^\ast P_{\vec{v}}}.
\]
So this is nothing other than the Gelfand-Naimark-Segal construction. Note that this inner product need not be positive; it will be positive (resp, positive definite) if the functional $\phi$ is. Define the action of $\mf{C} \langle X_1, X_2, \ldots, X_n \rangle$ on $\mf{C} \langle x_1, x_2, \ldots, x_n \rangle$ by \[
P_{\vec{u}} \cdot 1 = x_{\vec{u}}
\]
and
\[
P_{\vec{u}} \cdot x_{\vec{v}} = (P_{\vec{u}} P_{\vec{v}}) \cdot 1.
\]
Then for any $P$,
\[
\state{P} = \ip{1}{P \cdot 1}.
\]
$P_{\vec{u}}(\mb{X})$ is sometimes called the Wick product of $X_{u(1)}, X_{u(2)}, \ldots, X_{u(k)}$, and denoted $:X_{\vec{u}}:$.

\subsection{Appell Fock space}
Given a family of polynomials $\Aa{\vec{u}}{x_1, x_2, \ldots, x_n}$ satisfying
\[
\partial_{x_i} {A}_{\vec{u}} = \sum_{j: u(j) = i} {A}_{u(1), \ldots, \widehat{u(j)}, \ldots, u(k)},
\]
there is a unique functional $\phi$ on $\mf{R}[x_1, x_2, \ldots, x_n]$ such that $\state{1} = 1$, $\state{\Aa{\vec{u}}{\mb{x}}} = 0$ for all $\vec{u}$. The Fock space construction provides a natural way to recover such a functional. The induced inner product on $\mf{C} \langle x_1, x_2, \ldots, x_n \rangle$ is determined by
\[
\ip{x_{\vec{u}}}{x_{\vec{v}}} = \sum_{\substack{U \subset \set{1, \ldots, k}, U \neq \emptyset \\ V \subset \set{1, \ldots, l}, V \neq \emptyset}} R[X_{((\vec{u})^{op}:U)}, X_{(\vec{v}:V)}]
\]
for $\vec{u} \in \mf{N}^k$, $\vec{v} \in \mf{N}^l$. Since in this section the cumulant maps $R$ are symmetric in their arguments, the inner product is degenerate and factors through to $\mf{C}[x_1, x_2, \ldots, x_n]$.

\br
From the recursion relation~\eqref{Classical-recursion}, the action of the operator $X_j$ is
\[
X_j x_{\vec{u}} = x_j x_{\vec{u}} + \sum_{V \subset \set{1, \ldots, n}} R[X_j, X_{(\vec{u}:V)}] x_{(\vec{u}:V^c)}.
\]
Thus it is a sum of a creation operator
\[
a_j^\ast: x_{\vec{u}} \mapsto x_j x_{\vec{u}},
\]
a scalar operator
\[
x_{\vec{u}} \mapsto R[X_j] x_{\vec{u}},
\]
and some unusual annihilation operators.

\subsection{Sheffer and Meixner families}
Sheffer families are monic polynomial families $\set{P_n(x,t)}_{n=0}^\infty$ such that their exponential generating function has the form
\[
H(x,t,z) = \sum_{n=0}^\infty \frac{1}{n!} P_n(x,t) z^n = f(z)^t e^{U(z) x}
\]
for $U(z) = z + \text{ higher-order terms}$, $f(z) = 1 + \text{ higher-order terms}$. Among these, the Meixner families are those consisting of orthogonal polynomials (see \cite[Chapter 4]{SchOrthogonal} for references). It follows from the results of Meixner~\cite{Meixner} that up to affine transformations, the polynomials from this class satisfy recursion relations
\[
x P_n(x,t) = P_{n+1}(x,t) + a n P_n(x,t) + n (t + b (n-1)) P_{n-1}(x,t)
\]
for $a \in \mf{R}$, $b \in \mf{R}_+$. Here we have assumed, without loss of generality, that $P_1(x,t) = x$. By re-normalizing, we can restrict the values of the parameters to the following five classes, labelled by the names of the corresponding families.
\begin{description}
\item[Hermite] $a = b = 0$. Orthogonal with respect to the Gaussian distribution.
\item[Charlier] $b = 0$, $a = 1$. Orthogonal with respect to the centered Poisson distribution.
\item[Meixner-Pollaczek] $0 \leq a < 2$, $b=1$. Orthogonal with respect to the centered continuous binomial distribution.
\item[Laguerre] $a=2$, $b=1$. Orthogonal with respect to the centered Gamma distribution.
\item[Meixner] $a>2$, $b=1$. Orthogonal with respect to the centered negative binomial distribution.
\end{description}
Here for any measure $\mu$ with $m_1(\mu) < \infty$, the corresponding centered measure is
\[
\mu_c(S) = \mu(S + m_1(\mu)),
\]
for which
\[
\int_{\mf{R}} x \,d\mu_c(x) = m_1(\mu_c) = 0.
\]

\subsection{Kailath-Segall polynomials}
\label{Subsec:Kailath}
Basic Kailath-Segall polynomials $W^{(n)}$ have appeared in the paper~\cite{Kailath-Segall}. They were defined via certain stochastic integrals, but the authors also showed that they satisfy a recursion
\begin{equation}
\label{KS-formula}
W^{(n+1)} = \sum_{k=0}^n (-1)^k \frac{n!}{(n-k)!} x_{k+1} W^{(n-k)} - r_1 W^{(n)}
\end{equation}
(with a different notation and slightly different normalization). We define multivariate Kailath-Segall polynomials by recursion \eqref{Multi-KS}; then $W^{(n)} = W_{1, 1, \ldots, 1}$. Note that this is a natural definition in the compound Poisson (``de Finetti'') case; there is an analogous definition for the more general (``Kolmogorov'') case of functionals with finite variance, see \cite{AnsQLevy}. That paper also details the stochastic integral connection in Section~3.1.

\br
The following more general definition comes naturally from a Fock space construction, see Section~\ref{Section:KS-Fock-space}.

\begin{Defn}
Let $\mc{A}_0$ be a commutative complex star-algebra without identity, and $\Exp{\cdot}$ a star-linear functional on it. Denote by $\mc{A}_0^{sa}$ the self-adjoint elements of the algebra. Let $\mc{A}$ be the complex unital star-algebra generated by commuting symbols $\set{X(f): f \in \mc{A}_0^{sa}}$ (and $1$) subject to the linearity relations
\[
X(\alpha f + \beta g) = \alpha X(f) + \beta X(g).
\]
Equivalently, $\mc{A}$ is the symmetric tensor algebra of $\mc{A}_0$. The star-operation on it is determined by the requirement that all $X(f)$, $f \in \mc{A}_0^{sa}$ are self-adjoint. For such $f_i$, define the Kailath-Segall polynomials by $\KS{f} = X(f) - \Exp{f}$ and
\begin{equation*}
\begin{split}
\KS{f, f_1, f_2, \ldots, f_n}
& = X(f) \KS{f_1, f_2, \ldots, f_n}
- \sum_{i=1}^n \Exp{f f_i} \KS{f_1, \ldots, \hat{f}_i, \ldots, f_n} \\
&\quad - \sum_{i=1}^n \KS{f f_i, \ldots, \hat{f}_i, \ldots, f_n}
- \Exp{f} \KS{f_1, f_2, \ldots, f_n}.
\end{split}
\end{equation*}
The definition can be extended in a $\mf{C}$-linear way, so each $W$ is really a multi-linear map from $\mc{A}_0$ to $\mc{A}$, which turns out to be symmetric in its arguments.
\end{Defn}

\noindent
In the particular case $\mc{A}_0 = \mf{C}_0[x]$ (polynomials without constant term), we may denote $x_i = X(x^i)$. The functional can be taken to be the cumulant functional of a measure $\mu$, $\Exp{x^i} = r_i(\mu)$. Then
\[
W_{\vec{u}}(\mb{x}) = \KS{x^{u(1)}, x^{u(2)}, \ldots, x^{u(n)}}
\]
are multivariate polynomials in $\set{x_i: i \in \mf{N}}$. If $\phi$ is the functional on $\mf{C}[x_1, x_2, \ldots]$ determined by
\[
\state{W_{\vec{u}}(\mb{x})} = 0
\]
for all non-empty $\vec{u}$, then its cumulants are
\[
r_k(x_n) = r_{n k}(\mu),
\]
and more generally
\[
R[x_{\vec{u}}] = r_{\abs{\vec{u}}}(\mu).
\]
If $\mu$ is infinitely divisible, the functional $\Exp{\cdot}$ is positive. It follows that the functional $\phi$ is also positive, see Section~\ref{Section:KS-Fock-space}.

\begin{Prop}
\label{Prop:Kailath-Segall}
The following expansions hold.
\begin{enumerate}
\item
Of usual products in terms of the Kailath-Segall polynomials:
\begin{equation*}
X(f_1) X(f_2) \ldots X(f_n) =
\sum_{\pi \in \Part(n)} \sum_{S \subset \pi} \prod_{B \in S^c} \Exp{ f_B} \cdot \KS{f_C: C \in S}.
\end{equation*}
\item
Of the Kailath-Segall polynomials:
\[
\KS{f_1, f_2,, \ldots, f_n}
= \sum_{\pi \in \Part(n)} \sum_{S \subset \Sing(\pi)} (-1)^{n - \abs{\pi} + \abs{S}} \prod_{B \in S} \Exp{f_B} \prod_{C \in S^c} (\abs{C} - 1)! X(f_C).
\]
\item
Of products of Kailath-Segall polynomials: for $\vec{u} = (\vec{u}_1, \ldots, \vec{u}_k)$, 
\begin{multline*}
\prod_{i=1}^k \KS{f_{u_i(1)}, f_{u_i(2)}, \ldots, f_{u_i(s(i))}} \\
= \sum_{\substack{\pi \in \Part(N) \\ \pi \wedge \pi_{s(1), s(2), \ldots, s(k)} = \hat{0}}} \sum_{\substack{S \subset \pi, \\ \Sing(S^c) = \emptyset}} \prod_{B \in S^c} \Exp{f_{(\vec{u}:B)}} \KS{f_{(\vec{u}:C)}: C \in S}.
\end{multline*}
\end{enumerate}
\end{Prop}

\begin{proof}
We only consider part (c), see Theorem~\ref{Thm:Q-expansions} for the remaining proofs. We use induction on the length of the longest $\vec{u}_i$. If all of them have length $1$, the desired statement is
\[
\prod_{i=1}^k \KS{f_i}
= \sum_{\pi \in \Part(k)} \sum_{\substack{S \subset \pi, \\ \Sing(S^c) = \emptyset}} \prod_{B \in S^c} \Exp{f_B} W(f_C : C \in S).
\]
This is just the sum in part (a) except that $S^c$ is not allowed to contain any singletons, and the result follows from part (a) and the fact that $\KS{f_i} = X(f_i) - \Exp{f_i}$. Now let $\vec{u}_i$ be the longest multi-index. Suppose $\vec{u}_i(1) = j$, and denote $\vec{v} = (u_i(2), u_i(3), \ldots)$. Using the defining recursion relation,
\begin{equation*}
\begin{split}
&\KS{f_{\vec{u}_1}} \KS{f_{\vec{u}_2}} \ldots \KS{f_{\vec{u}_i}} \ldots \KS{f_{\vec{u}_k}} \\
&\qquad = \KS{f_{\vec{u}_1}} \ldots \biggl( X(f_j) \KS{f_{\vec{v}}} - \sum_{l=2}^{s(i)} \Exp{f_j f_{u_i(i)}} \KS{f_{u_i(1)}, \ldots, \hat{f}_{u_i(l)}, \ldots, f_{u_i(s(i))}} \\
&\qquad\quad - \sum_{l=2}^n \KS{f_j f_{u_i(l)}, \ldots, \hat{f}_{u_i(l)}, \ldots, f_{u_i(s(i))}} - \Exp{f_j} \KS{f_{\vec{v}}} \biggr) \ldots \KS{f_{\vec{u}_k}} \\
&\qquad = \KS{f_{\vec{u}_1}} \ldots \KS{f_j} \KS{f_{\vec{v}}} \ldots \KS{f_{\vec{u}_k}} \\
&\qquad\quad - \KS{f_{\vec{u}_1}} \ldots \sum_{l=2}^{s(i)} \Exp{f_j f_{u_i(l)}} \KS{f_{u_i(1)}, \ldots, \hat{f}_{u_i(l)}, \ldots, f_{u_i(s(i))}} \ldots \KS{f_{\vec{u}_k}} \\
&\qquad\quad - \KS{f_{\vec{u}_1}} \ldots \sum_{l=2}^n \KS{f_j f_{u_i(l)}, \ldots, \hat{f}_{u_i(l)}, \ldots, f_{u_i(s(i))}}  \ldots \KS{f_{\vec{u}_k}}
\end{split}
\end{equation*}
Apply the induction hypothesis to the right-hand-side. The desired sum on the left-hand-side is over pairs $(S, \pi)$, $\pi$ inhomogeneous. Any such term appears in the sum corresponding to the first term on the right-hand-side. We need to show that all the other elements in the sum corresponding to this term cancel out. Indeed, take any pair $S \subset \pi$, $\pi$ inhomogeneous with respect to $\pi_{s(1), \ldots, 1, s(i) -1, \ldots, s(k)}$, but not with respect to $\pi_{s(1), \ldots, s(i), \ldots, s(k)}$. Any such partition contains a class $B$ which contains the position of $u_i(1)$ and a position of some other $u_i(l)$. If $B \not \in S$ and $\abs{B} = 2$, it gets cancelled by the corresponding term from the second term on the right-hand-side. Otherwise it gets cancelled by the third term on the right-hand-side. 
\end{proof}

\begin{Cor}
\label{Cor:KS-linear}
The linearization coefficients for the Kailath-Segall polynomials are sums over inhomogeneous partitions with no singletons:
\[
\state{\prod_{i=1}^k \KS{f_{u_i(1)}, f_{u_i(2)}, \ldots, f_{u_i(s(i))}}}
= \sum_{\substack{\pi \in \Part(N) \\ \pi \wedge \pi_{s(1), s(2), \ldots, s(k)} = \hat{0}, \\
\Sing(\pi) = \emptyset}} \prod_{B \in \pi} \Exp{f_{(\vec{u}:B)}}.
\]
\end{Cor}

\begin{Cor}
$\set{W_{\vec{u}}}$ are pseudo-orthogonal.
\end{Cor}

\begin{proof}
For $\vec{u} \in \mf{N}^n$, $\vec{v} \in \mf{N}^k$, by the preceding corollary
\[
\state{W_{\vec{u}}^\ast W_{\vec{v}}}
= \delta_{n k} \sum_{\sigma \in \Sym(n)} \prod_{i=1}^n \Exp{x^{u(i)} x^{v(\sigma(i))}}. \qedhere
\]
\end{proof}

\noindent
Now assume that $\Exp{\cdot}$ is positive, $\Exp{x^i} = m_i(\nu)$ for some positive measure $\nu$. (In this case $\mu$ is infinitely divisible, in fact a compound Poisson measure, and $\nu$ is the L\'{e}vy measure for the convolution semigroup $\set{\mu_t}$.) Let $\set{p_i}_{i=1}^\infty$ be the orthogonal polynomials with respect to $\nu$, and $y_i = X(p_i)$. Equivalently, $\set{y_i}$ is the orthogonalization of $\set{x_i}$ with respect to the inner product $\ip{x_i}{x_j} = r_{i+j}(\mu) = m_{i+j}(\nu)$.

\begin{Cor}
\label{Cor:KS-orth}
By a linear change of variables, define polynomials
\[
W^\perp_{\vec{u}}(\mb{x}) = W_{\vec{u}}(\mb{y}).
\]
Then $W^\perp_{\vec{u}}$ are orthogonal, that is,
\[
\vec{u} \neq \vec{v} \Rightarrow \state{\left(W^\perp_{\vec{u}}(\mb{x})\right)^\ast W^\perp_{\vec{v}}(\mb{x})} = 0.
\]
\end{Cor}

\noindent
Free and even $q$-analogs of Corollary~\ref{Cor:KS-orth} hold, derived from appropriate modifications of Corollary~\ref{Cor:KS-linear}. These properties are related to the ``generalized chaos decomposition property'' of \cite{SchChaotic}.

\br
Appell polynomials are linear combinations of the Kailath-Segall polynomials of the same degree. Note that a priori, such a linear combination is a multivariate polynomial, but in this case it turns out to depend only on a single variable.

\begin{Prop}
\label{Prop:Appell-KS}
For $x = x_1$,
\[
A^{(n)}(x)
= \sum_{\substack{\pi \in \Part(n) \\ \pi = \set{B_1, B_2, \ldots, B_k}}} W_{\abs{B_1}, \abs{B_2}, \ldots, \abs{B_k}}(\mb{x})
= \sum_{k=1}^n \sum_{\vec{u} \in \Delta_{k,n}} a_{\vec{u}} W_{\vec{u}} (\mb{x}),
\]
where $a_{\vec{u}} = \abs{\set{\pi \in \Part(n), \pi = \set{B_1, B_2, \ldots, B_k}, u(i) = \abs{B_i}, i = 1, 2, \ldots, k}}$. Since in this case the variables commute and $W_{\vec{u}}$ is symmetric in its indices, also
\[
A^{(n)}(x) = \sum_{p \vdash n} \frac{n!}{(1!)^{p_1} \ldots (n!)^{p_n} p_1! \ldots p_n!} W_{p_1 \, 1's, p_2 \, 2's, \ldots, p_n \, n's}(\mb{x}),
\]
where $p \vdash n$ is a number partition of the number $n$, $\sum_{i=1}^n i p_i = n$.
\end{Prop}

\begin{proof}
The Fock representation of Section~\ref{Section:Fock} is clearly faithful. So it will suffice to show that in the Fock representation of the Kailath-Segall polynomials,
\[
A^{(n)}(X) \cdot 1
= \sum_{\substack{\pi \in \Part(n) \\ \pi = \set{B_1, B_2, \ldots, B_k}}} W_{\abs{B_1}, \abs{B_2}, \ldots, \abs{B_k}}(\mb{x}) \cdot 1
= \sum_{\pi \in \Part(n)} \prod_{i=1}^{\abs{\pi}} x_{\abs{B_i}}.
\]
The proof will proceed by induction. Using the recursion relation for the Appell polynomials and induction,
\begin{equation}
\label{Appell-KS}
A^{(n+1)}(X) \cdot 1
= \sum_{\pi \in \Part(n)} X \cdot \prod_{i=1}^{\abs{\pi}} x_{\abs{B_i}}
- \sum_{k=0}^n \binom{n}{k} r_{n+1-k} \sum_{\sigma \in \Part(k)} \prod_{i=1}^{\abs{\sigma}} x_{\abs{B_i}}.
\end{equation}
The action of the operator $X$ is determined by the recursion relation \eqref{Multi-KS} for the Kailath-Segall polynomials, 
\[
\begin{split}
X \cdot x_{\vec{u}}
& = (X_1 W_{\vec{u}}(\mb{X})) \cdot 1
= \Bigl[W_{1, u(1), u(2), \ldots, u(k)} \\
&\quad + \sum_{i=1}^k r_{1 + u(i)} W_{u(1), \ldots, \widehat{u(i)}, \ldots, u(k)} + \sum_{i=1}^k W_{1 + u(i), u(1), \ldots, \widehat{u(i)}, \ldots, u(k)} + r_1 W_{u(1), u(2), \ldots, u(k)} \Bigr] \cdot 1 \\
& = x_1 x_{\vec{u}} + \sum_{i=1}^k r_{1 + u(i)} x_{u(1)} \ldots \widehat{x_{u(i)}} \ldots x_{u(k)} + \sum_{i=1}^k x_{1 + u(i)} x_{u(1)} \ldots \widehat{x_{u(i)}} \ldots x_{u(k)} + r_1 x_{\vec{u}}.
\end{split}
\]
Therefore expression \eqref{Appell-KS} is
\begin{equation*}
\begin{split}
&\sum_{\pi \in \Part(n)} \biggl[ x_1 \prod_{i=1}^{\abs{\pi}} x_{\abs{B_i}}
+ \sum_{i=1}^{\abs{\pi}} x_{\abs{B_i} + 1} x_{\abs{B_1}} \ldots \hat{x}_{\abs{B_i}} \ldots x_{\abs{B_{\abs{\pi}}}} \\
&\qquad\qquad + \sum_{i=1}^{\abs{\pi}} r_{\abs{B_i} + 1} x_{\abs{B_1}} \ldots \hat{x}_{\abs{B_i}} \ldots x_{\abs{B_{\abs{\pi}}}}
+ r_1 \prod_{i=1}^{\abs{\pi}} x_{\abs{B_i}} \biggr] - \sum_{k=0}^n \binom{n}{k} r_{n+1-k} \sum_{\sigma \in \Part(k)} \prod_{i=1}^{\abs{\sigma}} x_{\abs{B_i}}.
\end{split}
\end{equation*}
The first term in the sum, as well as the fourth term, are indexed by all partitions of $(n+1)$ whose first class is a singleton. The second term, as well as the third term, are indexed by all partitions of $(n+1)$ whose first class is not a singleton. Finally, the term which is subtracted is indexed by all partitions of $(n+1)$, with the cumulant factor corresponding to the first class of such a partition, and the binomial factor accounting for the choice of the remaining $n-k$ elements of that class (other than the element $1$). Therefore the third and fourth terms in the brackets cancel exactly all the terms that are subtracted, and the first and second terms add up to the desired sum over $\Part(n+1)$.
\end{proof}

\begin{Ex}
Hermite polynomials are orthogonal, Appell, and Kailath-Segall (with $x_i = 0$ for $i>1$). Charlier polynomials are orthogonal, Sheffer, and Kailath-Segall (with all $x_i = x$). Jacobi polynomials are orthogonal, Bernoulli polynomials are Appell, Abel polynomials are Sheffer.
\end{Ex}

\newpage
\section{The free analogs}
\label{Section:Free}

\subsection{Notation}

\subsubsection{Non-commutative power series}
In this section we will use power series in non-commuting variables. Most theorems about formal power series remain valid in this context. In particular, series with non-zero constant term have inverses with respect to multiplication:
\[
(1 - \sum_{\vec{u}} a_{\vec{u}} z_{\vec{u}}) (1 + \sum_{\vec{v}} b_{\vec{v}} z_{\vec{v}}) = 1
\]
for
\begin{equation}
\label{Power-series}
b_{\vec{w}}
= \sum_{\substack{(\vec{u}, \vec{v}) = \vec{w}, \\ \vec{v} \neq \vec{w}}} a_{\vec{u}} b_{\vec{v}}
= \sum_{k=1}^n \sum_{(\vec{u}_1, \vec{u}_2, \ldots, \vec{u}_k) = \vec{w}} a_{\vec{u}_1} a_{\vec{u}_2} \ldots a_{\vec{u}_k},
\end{equation}
where $\vec{w} \in \mf{N}^n$. From general theory, left inverse and right inverse are equal.

\subsubsection{Free cumulants}
Let $G(z) = \sum_{n=0}^\infty m_n z^{-(n+1)}$ be a formal Laurent series, a generating function for a moment sequence. Define the corresponding free cumulant generating function $\mb{R}(z)$ via the functional relation
\begin{equation}
\label{Cauchy}
G \left(\frac{1 + \mb{R}(z)}{z} \right) = z; \qquad \frac{1 + \mb{R}(G(z))}{G(z)} = z.
\end{equation}
Note that we use the boldface notation to distinguish $\mb{R}$ from the usual $R$-transform $R(z) = \mb{R}(z)/z$ \cite{VoiSF}. Define the corresponding free cumulant sequence via its generating function
\[
\mb{R}(z) = \sum_{n=1}^\infty r_n z^n.
\]
Let $\mc{A}$ be a possibly non-commutative real algebra, and $\phi$ a unital real linear functional on it. For $\mb{X} = (X_1, X_2, \ldots, X_n) \subset \mc{A}$ a collection of non-commutative random variables, define their joint moments
\[
M[X_{\vec{u}}] = \state{X_{\vec{u}}}
\]
as before, but now
\[
\mb{M}(\mb{z})
= \sum_{\vec{u}} M[X_{\vec{u}}] z_{\vec{u}}
\]
denotes the \emph{ordinary} moment-generating function. Here, and in the sequel, $z_i$'s are formal non-commuting indeterminates. Define the joint free cumulants of $\mb{X}$ via
\begin{equation}
\label{Free-moment-cumulant}
R[X_{\vec{u}}]
= M[X_{\vec{u}}] - \sum_{\substack{\pi \in \NC(n), \\ \pi \neq \hat{1}}} \prod_{B \in \pi} R[X_{(\vec{u}:B)}],
\end{equation}
which expresses $R[X_{u(1)}, X_{u(2)}, \ldots, X_{u(n)}]$ in terms of the joint moments and sums of products of lower-order cumulants. Define the free cumulant generating function via
\[
\mb{R}(\mb{z})
= \sum_{\vec{u}} R[X_{\vec{u}}] z_{\vec{u}}.
\]
For a single random variable $X$, $R[X, X, \ldots, X] = r_n$ for the moment sequence $\set{m_k = \state{X^k}}$.

\br
The relation between joint moments and joint free cumulants can be summarized in a relation between their generating functions, as follows. The following proposition is due to Nica and Speicher; for completeness, we provide a direct proof.
\begin{Prop}\cite{Nica-Lecture}
\label{Prop:Nica}
Let $z_i = w_i \bigl(1 + \mb{M}(\mb{w}) \bigr)$. Then
\[
\mb{R}(\mb{z})
= \mb{R} \Bigl(w_1 \bigl(1 + \mb{M}(\mb{w}) \bigr), \ldots, w_n \bigl(1 + \mb{M}(\mb{w}) \bigr) \Bigr)
= \mb{M}(\mb{w}).
\]
\end{Prop}

\begin{proof}
In the defining relation~\eqref{Free-moment-cumulant}, the sum is over all non-crossing partitions. Each non-crossing partition can be described by a subset $V$ containing $1$ and a collection of non-crossing partitions on the intervals into which $V$ subdivides the set $\set{1, 2, \ldots, k}$. Applying the formula~\eqref{Free-moment-cumulant} again to each of those intervals, we obtain
\begin{equation*}
M[X_{\vec{u}}]
= \sum_{s=1}^\infty \sum_{\substack{V \subset \set{1, 2, \ldots, k} \\ V = \set{1 = i(1), i(2), \ldots, i(s)}}}
R[X_{(\vec{u}:V)}] \prod_{j=1}^s M[X_{\vec{u}:\set{i(j) + 1, \ldots, i(j+1) - 1}}],
\end{equation*}
where by convention $i(s+1) = k+1$ and $M[X_{\vec{u}:\set{i(j) + 1, \ldots, i(j+1) - 1}}] = M[\emptyset] = 1$ if $i(j+1) = i(j) + 1$. So
\begin{equation*}
\begin{split}
& \mb{M}(\mb{w})
= \sum_{\vec{u}} M[X_{\vec{u}}] w_{\vec{u}}
= \sum_{s=1}^\infty \sum_{\vec{v} \in \set{1, \ldots, n}^s} R[X_{\vec{v}}]
\sum_{\substack{\vec{v}_1, \vec{v}_2, \ldots, \vec{v}_k \\ \text{possibly empty}}}
\prod_{j=1}^s \left( w_{v(j)} M[X_{\vec{v}_j}] w_{\vec{v}_j} \right) \\
&\quad = \sum_{s=1}^\infty \sum_{\vec{v} \in \set{1, \ldots, n}^s} R[X_{\vec{v}}] \prod_{j=1}^s \left( w_{v(j)} (1 + \mb{M}(\mb{w})) \right)
= \mb{R} \Bigl(w_1 \bigl(1 + \mb{M}(\mb{w}) \bigr), \ldots, w_n \bigl(1 + \mb{M}(\mb{w}) \bigr) \Bigr).
\end{split}
\end{equation*}
\end{proof}

\begin{Ex}
The key distribution in free probability, which appears for example in the free central limit theorem, is the semicircular distribution. It is characterized by the property that $\mb{R}(\mb{z}) = z^2$, its moments are the Catalan numbers, and the corresponding measure has density
\[
\frac{1}{2 \pi} \sqrt{4 - x^2} \chf{[-2, 2]}(x) \,dx.
\]
More generally, we say that $X$ has a (scaled, shifted) semicircular distribution if $\mb{R}(\mb{z}) = a z + b z^2$. A family $\set{X_1, X_2, \ldots, X_n}$ form a semicircular family if all of them are self-adjoint and all their linear combinations have (scaled, shifted) semicircular distributions. Equivalently, they are self-adjoint and their free cumulant generating function is quadratic. They form a free semicircular family if in addition they are freely independent, in which case it suffices to require that each $X_i$ be  self-adjoint and have a semicircular distribution.
\end{Ex}

\subsection{Free Appell polynomials}
\begin{Defn}
\label{Defn:Free-Appell}
Free Appell polynomials are defined via their ordinary generating function
\[
\sum_{n=0}^\infty A^{(n)}(x) z^n
= \frac{1}{1 - x z + \mb{R}(z)}.
\]
\end{Defn}

\begin{Prop}
Some properties of the free Appell polynomials are:
\begin{gather}
\label{Free-Recursion}
x A^{(n)}(x) = A^{(n+1)} + \sum_{k=0}^n r_{n+1-k} A^{(k)}(x), \\
\notag
x^n
= \sum_{k=0}^n \biggl( \sum_{\vec{u} \in \Delta_{k+1, n+1}} \prod_{i=0}^k m_{u(i)-1} \biggr) A^{(k)}(x), \\
\notag
A^{(n)}(x) = \sum_{k=0}^n \biggl( 1 + \sum_{s=1}^{n-k} \sum_{\vec{u} \in \Delta_{s, n-k}} \binom{k+s}{s} (-1)^s r_{\vec{u}} \biggr) x^k.
\end{gather}
\end{Prop}

\br
See Proposition~\ref{Prop:Free-Appell} for a more general statement.

\begin{Ex}
Chebyshev polynomials of the second kind, $U_n(2 \cos \theta) = \frac{\sin (n+1) \theta}{\sin \theta}$, are free Appell. They are orthogonal with respect to the semicircular distribution. In fact, from equation~\eqref{Free-Recursion} it follows that these are the only orthogonal free Appell polynomials.
\end{Ex}

\subsection{Multivariate free Appell polynomials}
Define the partial derivative $\partial_i$ on polynomials in non-commuting variables $x_1, x_2, \ldots, x_n$ by a linear extension of its action on monomials
\[
\partial_i x_{\vec{u}} = \sum_{j: u(j) = i} x_{u(1)} \ldots \widehat{x_{u(j)}} \ldots x_{u(k)}.
\]

\begin{Defn}
\label{Defn:Multi-free-Appell}
Let $\mc{A}$ be a non-commutative real algebra with a unital real linear functional $\phi$.  We will continue to call its elements (non-commuting) random variables. For $X_1, X_2, \ldots, X_n \in \mc{A}$, define the multivariate free Appell polynomial
\[
\Aa{X_1, X_2, \ldots, X_n}{x_1, x_2, \ldots, x_n}
\]
by the following conditions:
\begin{equation}
\label{Appell-derivative}
\partial_i \Aa{X_1, X_2, \ldots, X_n}{x_1, \ldots, x_i, \ldots, x_n}
= \Aa{X_1, \ldots, X_{i-1}}{x_1, \ldots, x_{i-1}} \cdot \Aa{X_{i+1}, \ldots, X_n}{x_{i+1}, \ldots, x_n}
\end{equation}
(with $A_\emptyset = 1$),
\begin{equation}
\label{Expectation}
\state{\Aa{X_1, X_2, \ldots, X_n}{X_1, X_2, \ldots, X_n}} = 0,
\end{equation}
and the restriction that 
\begin{equation}
\label{Order}
\text{in each monomial, the variables appear in the increasing order of indices.}
\end{equation}
The polynomials are uniquely determined, and since for $i < j$,
\begin{equation*}
\begin{split}
&\partial_i \partial_j \Aa{X_1, X_2, \ldots, X_n}{x_1, \ldots, x_i, \ldots, x_n} \\
&\quad = \Aa{X_1, \ldots, X_{i-1}}{x_1, \ldots, x_{i-1}} \cdot \Aa{X_{i+1}, \ldots, X_{j-1}}{x_{i+1}, \ldots, x_{j-1}} \cdot \Aa{X_{j+1}, \ldots, X_n}{x_{j+1}, \ldots, x_n} \\
&\quad = \partial_j \partial_i \Aa{X_1, X_2, \ldots, X_n}{x_1, \ldots, x_i, \ldots, x_n},
\end{split}
\end{equation*}
they are well defined. Most of the time we will be interested in
\[
\Aa{X_1, X_2, \ldots, X_n}{X_1, X_2, \ldots, X_n},
\]
which will be denoted simply by $\A{X_1, X_2, \ldots, X_n}$.

\br
On the other hand, if the $n$-tuple $\set{X_1, X_2, \ldots, X_n}$ is fixed, define the polynomials $\Aa{\vec{u}}{\mb{x}}$ by
\[
\Aa{\vec{u}}{x_i: i = u(j) \text{ for some } j} = \Aa{X_{u(1)}, X_{u(2)}, \ldots, X_{u(k)}}{x_{u(1)}, x_{u(2)}, \ldots, x_{u(k)}}.
\]
This is a monic polynomial family.
\end{Defn}

\begin{Ex}
From conditions~\eqref{Appell-derivative} and~\eqref{Expectation} is follows that each $\A{X_1, \ldots, X_n}$ is a polynomial in $X_1, \ldots, X_n$ and their joint moments. Moreover, from conditions~\eqref{Appell-derivative} and~\eqref{Order} it follows that this polynomial has degree $1$ in each $X_i$. Some low-order polynomials are:
\begin{align*}
\A{X_1} &= X_1 - \state{X_1}, \\
\A{X_1, X_2} &= X_1 X_2 - \state{X_1} X_2 - \state{X_2} X_1 - \state{X_1 X_2} + 2 \state{X_1} \state{X_2}, \\
\A{X_1, X_2, X_3}
&= X_1 X_2 X_3 - \state{X_1} X_2 X_3 - \state{X_2} X_1 X_3 - \state{X_3} X_1 X_2 \\
&\quad - \state{X_1 X_2} X_3 - \state{X_2 X_3} X_1 + 2 \state{X_1} \state{X_2} X_3 + \state{X_1} \state{X_3} X_2 \\
&\quad + 2 \state{X_2} \state{X_3} X_1 - \state{X_1 X_2 X_3} + 2 \state{X_1 X_2} \state{X_3} + 2 \state{X_1} \state{X_2 X_3} \\
&\quad+ \state{X_1 X_3} \state{X_2} - 5 \state{X_1} \state{X_2} \state{X_3}.
\end{align*}
\end{Ex}

\begin{Lemma}
For $n \geq 1$, the map $\mc{A}^n \rightarrow \mc{A}$,
\[
(X_1, X_2, \ldots, X_n) \mapsto \A{X_1, X_2, \ldots, X_n}
\]
is $n$-linear.
\end{Lemma}

\begin{proof}
For $X_1 = \alpha X + \beta Y$, we will show that
\begin{multline}
\label{Linearity}
\Aa{\alpha X + \beta Y, X_2, \ldots, X_n}{\alpha x + \beta y, x_2, \ldots, x_n} \\
= \alpha \Aa{X, X_2, \ldots, X_n}{x, x_2, \ldots, x_n} + \beta \Aa{Y, X_2, \ldots, X_n}{y, x_2, \ldots, x_n}.
\end{multline}
Then in particular,
\[
\A{\alpha X + \beta Y, X_2, \ldots, X_n}
= \alpha \A{X, X_2, \ldots, X_n} + \beta \A{Y, X_2, \ldots, X_n},
\]
and a similar proof holds for the other components.
\begin{multline*}
\partial_x \Aa{\alpha X + \beta Y, X_2, \ldots, X_n}{\alpha x + \beta y, x_2, \ldots, x_i, \ldots, x_n}
= \alpha \Aa{X_2, \ldots, X_n}{x_2, \ldots, x_n} \\
= \partial_x \bigl( \alpha \Aa{X, X_2, \ldots, X_n}{x, x_2, \ldots, x_n}
+ \beta \Aa{Y, X_2, \ldots, X_n}{y, x_2, \ldots, x_n} \bigr).
\end{multline*}
So condition~\eqref{Order} implies that each monomial containing $x$ in $\alpha \Aa{X, X_2, \ldots, X_n}{x, x_2, \ldots, x_n}$ appears, with the same coefficient, in $\Aa{\alpha X + \beta Y, X_2, \ldots, X_n}{\alpha x + \beta y, x_2, \ldots, x_i, \ldots, x_n}$, with $x$ replaced by $(\alpha x + \beta y)$, and all monomials containing $(\alpha x + \beta y)$ appear in this way. Application of $\partial_y$ produces a similar statement, and implies in particular that the right-hand-side of~\eqref{Linearity} is a polynomial in $(\alpha x + \beta y)$, and that monomials containing this term on both sides of ~\eqref{Linearity} coincide. By induction,
\[
\begin{split}
&\partial_i \Aa{\alpha X + \beta Y, X_2, \ldots, X_n}{\alpha x + \beta y, x_2, \ldots, x_i, \ldots, x_n} \\
&\quad= \Aa{\alpha X + \beta Y, \ldots, X_{i-1}}{\alpha x + \beta y, \ldots, x_{i-1}}
\cdot \Aa{X_{i+1}, \ldots, X_n}{x_{i+1}, \ldots, x_n} \\
&\quad= \alpha \Aa{X, \ldots, X_{i-1}}{x, \ldots, x_{i-1}}
\cdot \Aa{X_{i+1}, \ldots, X_n}{x_{i+1}, \ldots, x_n} \\
&\quad\quad + \beta \Aa{Y, \ldots, X_{i-1}}{y, \ldots, x_{i-1}}
\cdot \Aa{X_{i+1}, \ldots, X_n}{x_{i+1}, \ldots, x_n} \\
&\quad= \partial_i \Bigl( \alpha \Aa{X, X_2, \ldots, X_n}{x, x_2, \ldots, x_n}
+ \beta \Aa{Y, X_2, \ldots, X_n}{y, x_2, \ldots, x_n} \Bigr).
\end{split}
\]
This implies that all the terms on both sides of~\eqref{Linearity} which contain $\alpha x + \beta y, x_2, \ldots, x_n$ coincide. Finally,
\begin{multline*}
\state{\Aa{\alpha X + \beta Y, X_2, \ldots, X_n}{\alpha X + \beta Y, X_2, \ldots, X_n}} = 0 \\
= \state{\alpha \Aa{X, X_2, \ldots, X_n}{X, X_2, \ldots, X_n} + \beta \Aa{Y, X_2, \ldots, X_n}{X, X_2, \ldots, X_n}}.
\end{multline*}
implies that the constant terms coincide as well.
\end{proof}

\br
In the following power series, the $x$ variables commute with the $z$ variables, although neither the $x$ nor the $z$ variables commute among themselves.

\begin{Lemma}
\label{Lemma:B-R}
Let $\mb{X} = (X_1, X_2, \ldots, X_n)$ be a family of random variables, and $\mb{R}$ the corresponding free cumulant functional.  Then in the expansion
\[
\left( 1 - \mb{x} \cdot \mb{z} + \mb{R}(\mb{z}) \right)^{-1}
= 1 + \sum_{\vec{u}} B_{\vec{u}}(\mb{x}) z_{\vec{u}},
\]
\[
B_{\vec{u}}(\mb{x})
= \sum_{V \subset \set{1, \ldots, n}} x_{(\vec{u}:V)} \sum_{\substack{\pi \in \Int(V^c) \\ (\pi, V) \in \NC(n)}} (-1)^{\abs{\pi}} R_\pi[X_{(\vec{u}:V^c)}].
\]
\end{Lemma}

\begin{proof}
This is just formula \eqref{Power-series} for the coefficients of the inverse power series to
\[
1 - \sum_{i=1}^n x_i z_i - \sum_{\vec{u}} \bigl(- R[X_{\vec{u}}]\bigr) z_{\vec{u}}.
\]
\end{proof}

\begin{Thm}
Let $\mb{X} = (X_1, X_2, \ldots, X_n)$ be a family of random variables, and $\mb{M}, \mb{R}, A$ be the corresponding moment and free cumulant functionals and free Appell polynomials. Then
\[
1 + \sum_{\vec{u}} \Aa{\vec{u}}{\mb{x}} z_{\vec{u}}
= H(\mb{x}, \mb{z})
= \left( 1 - \mb{x} \cdot \mb{z} + \mb{R}(\mb{z}) \right)^{-1}.
\]
\end{Thm}

\begin{proof}
Using the substitution in Proposition~\ref{Prop:Nica},
\begin{equation*}
\begin{split}
H(\mb{x}, \mb{z})
&= \biggl( 1 + R \Bigl(w_1 \bigl(1 + \mb{M}(\mb{w}) \bigr), \ldots, w_n \bigl(1 + \mb{M}(\mb{w}) \bigr) \Bigr) - \mb{x} \cdot \mb{w} \bigl(1 + \mb{M}(\mb{w}) \bigr) \biggr)^{-1} \\
&= \Bigl( 1 + \mb{M}(\mb{w}) - \mb{x} \cdot \mb{w} \bigl(1 + \mb{M}(\mb{w}) \bigr) \Bigr)^{-1} \\
&= (1 + \mb{M}(\mb{w}))^{-1} (1 - \mb{x} \cdot \mb{w})^{-1}.
\end{split}
\end{equation*}
Since
\[
(1 - \mb{x} \cdot \mb{w})^{-1} = 1 + \sum_{\vec{u}} x_{\vec{u}} w_{\vec{u}},
\]
it follows that
\[
\state{H(\mb{X}, \mb{z})}
= (1 + \mb{M}(\mb{w}))^{-1} \state{(1 - \mb{X} \cdot \mb{w})^{-1}}
= 1.
\]
So in the notation of Lemma~\ref{Lemma:B-R}, $\state{B_{\vec{u}}(\mb{X})} = 0$.

\br
Since
\[
\partial_{x_i} (1 - \mb{x} \cdot \mb{z} + \mb{R}(\mb{z}))^{-1}
= (1 - \mb{x} \cdot \mb{z} + \mb{R}(\mb{z}))^{-1} z_i (1 - \mb{x} \cdot \mb{z} + \mb{R}(\mb{z}))^{-1},
\]
\[
\partial_{x_i} {B}_{\vec{u}} = \sum_{j: u(j) = i} {B}_{u(1), \ldots, u(j-1)} {B}_{u(j+1), \ldots, u(k)}.
\]

\br
Finally, it follows from Lemma~\ref{Lemma:B-R} that in each monomial of $B_{1,2, \ldots, n}$, the variables appear in the increasing order of the indices. Therefore by the definition of free Appell polynomials, $B_{1,2, \ldots, n} = A_{1,2, \ldots, n}$.

\br
This is true for an arbitrary family of random variables, in particular for
\[
(Y_1, Y_2, \ldots, Y_k) 
= (X_{u(1)}, X_{u(2)}, \ldots, X_{u(k)}).
\]
Using the explicit expression from Lemma~\ref{Lemma:B-R}, $B_{\vec{u}}(\mb{x})$ is the free Appell polynomial $A_{1,2,\ldots, k}$ for $\mb{Y}$ applied to $(x_{u(1)}, \ldots, x_{u(k)})$, which is exactly $A_{\vec{u}}$ (for $\mb{X}$).
\end{proof}

\begin{Cor}
The one-variable polynomials of Definition~\ref{Defn:Free-Appell} are $A^{(n)}(x_1) = A_{1,1, \ldots, 1}(x_1)$. 
\end{Cor}

\begin{Prop}
\label{Prop:Free-Appell}
Let $\vec{u} \in \mf{N}^n$.
\begin{enumerate}
\item
The expansion of $X_{\vec{u}}$ in terms of the free Appell polynomials is
\begin{equation}
\label{Free-X-Appell}
\begin{split}
X_{\vec{u}}
&= \sum_{k=0}^n \sum_{\substack{B \subset \set{1, 2, \ldots, n}, \\ B =
\set{v(1) < v(2) < \ldots < v(k)}}} \A{X_{(\vec{u}:B)}}
\prod_{j=1}^{k+1}M[X_{(\vec{u}:\set{v(j-1) + 1, \ldots, v(j) - 1})}] \\
&= \sum_{\pi \in \NC(n)} \biggl( \prod_{C \in \pi} R[X_{(\vec{u}:C)}] +
\sum_{B \in \Outer(\pi)} \A{X_{(\vec{u}:B)}} \prod_{\substack{C \in \pi,
\\ C \neq B}} R[X_{(\vec{u}:C)}] \biggr),
\end{split}
\end{equation}
where $v(0) = 0$, $v(k+1) = n+1$,
\item
The recursion relation for the free Appell polynomials is
\begin{equation}
\label{Free-Appell-recursion}
\A{X_{(j, \vec{u})}} = X_j \A{X_{\vec{u}}} - \sum_{k=0}^n R[X_j, X_{u(1)}, \ldots, X_{u(n-k)}] \A{X_{u(n-k+1)}, \ldots, X_{u(n)}}.
\end{equation}
\item
The expansion of the free Appell polynomial is
\[
\A{X_{\vec{u}}} = \sum_{V \subset \set{1, \ldots, n}} X_{(\vec{u}:V)} \sum_{\substack{\pi \in \Int(V^c) \\ (\pi, V) \in \NC(n)}} (-1)^{\abs{\pi}} R_\pi[X_{(\vec{u}:V^c)}].
\]
\end{enumerate}
\end{Prop}

\begin{proof}
From the proof of the preceding theorem,
\[
(1 - \mb{x} \cdot \mb{w})^{-1}
= \bigl(1 + \mb{M}(\mb{w}) \bigr) (1 + \sum_{\vec{u}} \Aa{\vec{u}}{\mb{x}} z_{\vec{u}}),
\]
where
\[
z_i = w_i \bigl(1 + \mb{M}(\mb{w}) \bigr).
\]
But
\[
(1 - \mb{x} \cdot \mb{w})^{-1} = 1 + \sum_{\vec{u}} x_{\vec{u}} w_{\vec{u}}.
\]
So $x_{\vec{u}}$ is the coefficient of $w_{\vec{u}}$ in
\[
\sum_{k=0}^\infty \sum_{\vec{v} \in \mf{N}^k} \Aa{\vec{v}}{\mb{x}} \bigl(1
+ \mb{M}(\mb{w}) \bigr) w_{v(1)} \bigl(1 + \mb{M}(\mb{w}) \bigr) \ldots
w_{v(k)} \bigl(1 + \mb{M}(\mb{w}) \bigr).
\]
The first line of equation~\eqref{Free-X-Appell} follows. The second line, including the fact that $B$ is outer, follows from the definition of the free cumulants.

\br
The recursion relation follows by expanding the identity
\[
(1 - \mb{x} \cdot \mb{z} + \mb{R}(\mb{z})) (1 + \sum_{\vec{u}} \Aa{\vec{u}}{\mb{x}} z_{\vec{u}})
= 1.
\]
The final expansion follows by combining the preceding theorem with Lemma~\ref{Lemma:B-R}.
\end{proof}

\subsection{Free binomial}
Suppose $\pi \in \Part(n)$ has the property that the collections $\set{X_i: i \in B}_{B \in \pi}$ are freely independent. Then the basic relation between free independence and free cumulants \cite{SpeNCReview} says that
\begin{equation}
\label{Free-sum}
\mb{R}(\mb{z}) = \sum_{B \in \pi} \mb{R}(z_i: i \in B).
\end{equation}
More precisely, if $\mb{R}$ is the free cumulant generating function for $\set{X_1, \ldots, X_n}$ (and so a function of $z_1, \ldots, z_n$), and for each subset $B \subset \set{1, \ldots, n}$, $\mb{R}_B$ is the free cumulant generating function for $\set{X_i: i \in B}$ (and so a function of $\set{z_i: i \in B}$), then $\mb{R} = \sum_{B \in \pi} \mb{R}_B$.

\br
Now let $\vec{u} \in \Delta_{k, n}$. It can be written uniquely as $\vec{u} = (\vec{v}_1, \vec{v}_2, \ldots, \vec{v}_k)$ with each $\vec{v}_i$ longest consecutive sequence from the same class of $\pi$.
\begin{Prop}
\label{Prop:binomial}
In the setting above,
\[
\A{X_{\vec{u}}} = \prod_{i=1}^k \A{X_{\vec{v}_i}}.
\]
\end{Prop}

\noindent
The proof is based on the following lemma due to Franz Lehner. The lemma deals with bounded operators, but it applies equally well to formal power series.

\begin{Lemma} \cite{LehSpectra}
\label{Lemma:Lehner}
Let $S_1, \ldots, S_N \in B(H)$ be arbitrary operators and assume that the sum of alternating products
\[
I + \sum_{n=1}^\infty \sum_{i_1 \neq i_2 \neq \ldots \neq i_n} S_{i_1} S_{i_2} \ldots S_{i_n}
\]
(the sum over all products where neighboring factors are different) converges, then it equals
\[
\left( I - \sum_{i=1}^N S_i(I + S_i)^{-1} \right)^{-1}.
\]
\end{Lemma}

\begin{proof}[Proof of Proposition~\ref{Prop:binomial}]
Let $\pi$ have classes $B_1, B_2, \ldots, B_N$. From~\eqref{Free-sum},
\begin{equation}
\label{sum-R}
1 - \mb{X} \cdot \mb{z} + \mb{R}(\mb{z}) = 1 - \sum_{j=1}^N \Bigl[ \sum_{i \in B_j} X_i z_i - \mb{R}_{B_j} \Bigr].
\end{equation}
Denote
\[
H = \bigl( 1 - \mb{X} \cdot \mb{z} + \mb{R}(\mb{z}) \bigr)^{-1} = 1 + \sum_{k=1}^\infty \sum_{\vec{u} \in \set{1, \ldots, n}^k} \A{X_{\vec{u}}} z_{\vec{u}}
\]
and for $j=1, \ldots, N$,
\[
H_j = \Bigl( 1 - \sum_{i \in B_j} \mb{X} \cdot \mb{z} + \mb{R}_{B_j} \Bigr)^{-1} = 1 + \sum_{k=1}^\infty \sum_{\vec{v} \in (B_j)^k} \A{X_{\vec{v}}} z_{\vec{v}}.
\]
Using Lemma~\ref{Lemma:Lehner} with $S_j = H_j - 1$ and equation~\eqref{sum-R}, we get
\[
\begin{split}
1 + \sum_{n=1}^\infty \sum_{i_1 \neq i_2 \neq \ldots \neq i_n} (H_{i_1}-1) (H_{i_2}-1) \ldots (H_{i_n}-1)
& = \left( 1 - \sum_{i=1}^N (H_i-1) H_i^{-1} \right)^{-1} \\
& = \left( 1 - \sum_{i=1}^N (1 - H_i^{-1}) \right)^{-1} = H.
\end{split}
\]
Equating coefficients of $z_{\vec{u}}$ on both sides of this equation gives the proposition.
\end{proof}

\noindent
As a consequence, we get the ``non-commutative binomial'' formula: for $X, Y$ freely independent,
\[
\begin{split}
A^{(n)}(X + Y)
& = \A{X + Y, \ldots, X + Y} \\
& = \sum_{k=1}^n \sum_{\vec{u} \in \Delta_{k,n}} \Bigl[ A^{(u(1))}(X) A^{(u(2))}(Y) A^{(u(3))}(X)\ldots + A^{(u(1))}(Y) A^{(u(2))}(X) \ldots \Bigr].
\end{split}
\]
This is the extension to the free Appell polynomials of the binomial expansion of $(X + Y)^n$ for $X, Y$ non-commuting variables.

\br
A different analog of the commutative binomial expansion \eqref{binom} is the ``co-multiplication'' property of the free Appell polynomials: it easily follows from Definition~\ref{Defn:Free-Appell} that for general $x, y$,
\[
\frac{A^{(n)}(x) - A^{(n)}(y)}{x-y} = \sum_{k=0}^{n-1} A^{(k)}(x) A^{(n-k-1)}(y),
\]
where the $x$ and $y$-dependent polynomials correspond to possibly different functionals. In particular,
\[
\partial_x A^{(n)}(x) = \sum_{k=0}^{n-1} A^{(k)}(x) A^{(n-k-1)}(x).
\]

\subsection{Martingale property}
The martingale property of free Appell (and, more generally, Sheffer) polynomials for processes with freely independent increments was shown previously by Biane~\cite{BiaProcesses}. The following is an alternative proof for distributions all of whose moments are finite, which uses the binomial property above.
\begin{Lemma}
Let $\set{\mu_t}$ be a free convolution semigroup with all moments finite, $\set{X(t)}$ the corresponding free L\'{e}vy process, and $\set{A^{(n)}}$ the corresponding free Appell polynomials. Then for each $n$, the process $\set{A^{(n)}(X(t))}$ is a martingale with respect to the filtration induced by $\set{X(t)}$.
\end{Lemma}

\noindent
We again note that $A^{(n)}(X(t)) = A^{(n)}_{X(t)}(X(t))$, and the polynomial $A^{(n)}_{X(t)}$ depends on $t$, so the short-hand notation should be handled with care.

\begin{proof}
\[
\begin{split}
\Cexp{s}{A^{(n)}(X(t))}
&= \Cexp{s}{A^{(n)} \Bigl(X(s) + \bigl(X(t) - X(s) \bigr) \Bigr)} \\
&= \Cexp{s}{\sum_{k=1}^n \sum_{\vec{u} \in \Delta_{k,n}} A^{(u(1))} \bigl(X(s) \bigr) A^{(u(2))} \bigl(X(t) - X(s) \bigr) \ldots}
= A^{(n)}(X(s))
\end{split}
\]
since $\E \left[A^{(k)} \bigl(X(t) - X(s) \bigr) \right] = E \left[A^{(k)} \bigl(X(s) \bigr) \right] = 0$ for $k > 0$, and using the definition of free independence and properties of conditional expectation.
\end{proof}

\begin{Thm}
\label{Free-Appell-properties}
Let $\vec{u}_i \in \mf{N}^{s(i)}$, $i = 1, 2, \ldots, k$, $N = \sum_{i=1}^k s(i)$, and $\vec{u} = (\vec{u}_1, \vec{u}_2, \ldots, \vec{u}_k)$. Each of the quantities
\begin{gather}
\label{Moment-Full}
M[X_{\vec{u}_1}, X_{\vec{u}_2}, \ldots, X_{\vec{u}_k}]
= M[X_{\vec{u}}] \\
\label{Moment-Appell}
M \bigl[ \A{X_{\vec{u}_1}}, \A{X_{\vec{u}_2}}, \ldots, \A{X_{\vec{u}_k}} \bigr] \\
\label{Cumulant-Full}
R[(X_{\vec{u}_1}), (X_{\vec{u}_2}), \ldots, (X_{\vec{u}_k})] \\
\label{Cumulant-Appell}
R \bigl[ \A{X_{\vec{u}_1}}, \A{X_{\vec{u}_2}}, \ldots, \A{X_{\vec{u}_k}} \bigr]
\end{gather}
is equal to the sum of $R_{\pi}[X_{\vec{u}}]$ over partitions $\pi \in \NC(N)$ which are:
\begin{enumerate}
\item
Equation~\eqref{Moment-Full}: arbitrary,
\item
Equation~\eqref{Moment-Appell}: non-homogeneous,
\item
Equation~\eqref{Cumulant-Full}: connected,
\item
Equation~\eqref{Cumulant-Appell}: connected and non-homogeneous
\end{enumerate}
with respect to the interval partition $\pi_{s(1), s(2), \ldots, s(k)}$. Note that~\eqref{Moment-Appell} is the linearization coefficient for the free Appell polynomials.
\end{Thm}

\noindent
There is a similar theorem for the usual Appell polynomials, obtained by replacing the lattice of non-crossing partitions with the lattice of all partitions. The following proof is also very similar to the one in \cite{Giraitis-Surgailis}.

\begin{proof}
Equation~\eqref{Moment-Full} is the basic relation between moments and free cumulants. Equation~\eqref{Cumulant-Full} was proven in \cite{KraSpeProducts,Spe-Conceptual}. We prove equation~\eqref{Cumulant-Appell}; the proof of \eqref{Moment-Appell} is similar, and also implied by equation~\eqref{Appell-linearization} below. In fact, we will prove a more general statement, that
\begin{equation}
\label{Cumulant-General}
R \bigl[ \A{X_{\vec{u}_1}}, \ldots, \A{X_{\vec{u}_j}}, (X_{\vec{u}_{j+1}}), \ldots, (X_{\vec{u}_k}) \bigr]
\end{equation}
is equal to the sum of $R_{\pi}[X_{\vec{u}}]$ over all non-crossing connected partitions with no homogeneous classes in $\pi_{s(1), s(2), \ldots, s(j)}$. The proof will proceed by induction on $s(1), \ldots, s(j)$, starting with the statement for
\[
R \bigl[ X_{v(1)}, \ldots, X_{v(j)}, (X_{\vec{u}_{j+1}}), \ldots, (X_{\vec{u}_k}) \bigr]
\]
which is valid by equation~\eqref{Cumulant-Full}. Suppose \eqref{Cumulant-General} holds for all smaller $j$ or, for the same $j$, for all shorter $\vec{u}_1, \ldots, \vec{u}_j$. Substituting the recursion relation
\[
\A{X_{\vec{u}_j}}
= X_{\vec{u}_j} - \sum_{\substack{\pi \in \NC(s(j)) \\ \pi \neq \hat{1}}} \sum_{B \in \pi} \A{X_{(\vec{u}_j:B)}} \prod_{\substack{C \in \pi, \\ C \neq B}} R[X_{(\vec{u}_j:C)}]
\]
into equation~\eqref{Cumulant-General} breaks it up into the difference of two terms. The first term contains only $(j-1)$ Appell entries. It is the sum over all non-crossing connected partitions with no homogeneous classes in $\pi_{s(1), s(2), \ldots, s(j-1)}$. The second term contains $j$ Appell entries, the first $(j-1)$ of which have indices $\vec{u}_1, \ldots, \vec{u}_{j-1}$, and the last one is a proper sub-index of $\vec{u}_j$. For each choice of a proper subset 
\[
B \subset \set{\sum_{i=1}^{j-1} s(i) + 1, \ldots, \sum_{i=1}^j s(i)},
\]
the second term contains the sum over all non-crossing connected partitions such that the complement of $B$ in this subset consists of homogeneous classes, and $\pi$ has no homogeneous classes with respect to $\pi_{s(1), s(2), \ldots, B}$. In other words, the second term is the sum over all non-crossing connected partitions which have no homogeneous classes with respect to $\pi_{s(1), s(2), \ldots, s(j-1)}$ but have some homogeneous classes in the $j$'th class of $\pi_{s(1), s(2), \ldots, s(k)}$. Clearly the difference of these terms is the sum over all non-crossing connected partitions with no homogeneous classes in $\pi_{s(1), s(2), \ldots, s(j)}$.
\end{proof}

\noindent
The following proposition is a generalization of~\eqref{Moment-Appell}. Its proof is similar.

\begin{Prop}
Products of free Appell polynomials have an expansion
\begin{equation}
\label{Appell-linearization}
\A{X_{\vec{u}_1}} \A{X_{\vec{u}_2}} \ldots \A{X_{\vec{u}_k}}
= \sum_{V \subset \set{1, \ldots, N}} \sum_{\substack{\pi \in \NC(V^c) \\ \text{non-homogeneous}, \\ (\pi, V) \in \NC(N), \\ V \in \text{ Outer}(\pi, V)}} R_{\pi}[X_{(\vec{u}:V^c)}] \A{X_{(\vec{u}:V)}}.
\end{equation}
Here, ``non-homogeneous'' is with respect to the restriction $\pi_{s(1), s(2), \ldots, s(k)} \upharpoonright V^c$.
\end{Prop}

\begin{Prop}
Fix $\set{X_1, X_2, \ldots, X_n}$ with the moment generating function $\mb{M}$, free cumulant generating function $\mb{R}$, and free Appell polynomials $\set{{A}_{\vec{u}}}$. With respect to the joint distribution functional $\phi_{\mb{X}}$,
\begin{enumerate}
\item
$\set{{A}_{\vec{u}}}$ are pseudo-orthogonal iff $\set{X_1, X_2, \ldots, X_n}$ form a semicircular family.
\item
$\set{{A}_{\vec{u}}}$ are orthogonal iff $\set{X_1, X_2, \ldots, X_n}$ form a free semicircular family.
\end{enumerate}
\end{Prop}

\begin{proof}
Suppose $\set{A_{\vec{u}}}$ are pseudo-orthogonal. By formula~\eqref{Moment-Appell}, for $k>2$,
\[
\state{\A{X_{u(1)}}^\ast \A{X_{u(2)}, \ldots, X_{u(k)}}}
= R[X_{u(1)}, X_{u(2)}, \ldots, X_{u(k)}] = 0.
\]
So all the joint cumulants of $\set{X_1, X_2, \ldots, X_n}$ of order greater than $2$ are zero, and these random variables form a semicircular family. In particular, in formula~\eqref{Moment-Appell} only pair partitions make a non-zero contribution. If in addition $\set{A_{\vec{u}}}$ are orthogonal, from this formula it follows in addition that $R[X_i, X_j] = 0$ unless $i=j$. Therefore
\[
\mb{R}(\mb{z})
= \sum_{i=1}^n (R[X_i, X_i] z_i^2 + R[X_i] z_i).
\]
So $\set{X_1, X_2, \ldots, X_n}$ are freely independent, and form a free semicircular family.
\end{proof}

\begin{Remark}
I thank Dima Shlyakhtenko for bringing to my attention the following observation. Given a family of self-adjoint elements $\set{X_k}$ whose free Appell polynomials are pseudo-or\-tho\-go\-nal, there need not exist a linear change of variables $B$ such that $\set{Y_i = \sum_{j=1}^n B_{ij} X_j}$ is a family of self-adjoint elements whose free Appell polynomials are orthogonal. Indeed, one can always find such a linear transformation with $R[Y_i, Y_j] = \delta_{i j}$, but the $Y_i$'s need not be self-adjoint. As a result, $\state{Y_i^\ast Y_j}$ need not be $0$. In other words, while the joint free cumulants of the $Y_i$'s are $0$, the joint free star-cumulants of the $Y_i$'s need not be $0$, so the $Y_i$'s need not be freely independent.

\br
If $\phi$ has the trace property $\state{a b} = \state{b a}$, then one \emph{can} always orthogonalize pseudo-orthogonal polynomials.
\end{Remark}

\subsection{Free Sheffer and Meixner families}
Sheffer families are martingale polynomials for the corresponding L\'{e}vy processes. Based on this idea and the result of \cite{BiaProcesses}, in \cite{AnsMeixner} we proposed the definition of free Sheffer families, which are families of martingale polynomials with respect to free L\'{e}vy processes. Specifically, free Sheffer families are families of monic polynomials whose ordinary generating function $H(x,t,z) = \sum_{n=0}^\infty P_n(x,t) z^n$ has the form
\[
\frac{1}{1 + f(z) t - U(z) x} = \frac{1}{1 + f(z) t} \cdot \frac{1}{1 - \frac{U(z)}{1 + f(z) t} x}.
\]
We also described all the free Meixner families, that is, free Sheffer families consisting of orthogonal polynomials. They are given by the recursion relations
\begin{align*}
P_1(x,t) &= x, \\
x P_1(x,t) &= P_2(x,t) + a x + t, \\
x P_n(x,t) &= P_{n+1}(x,t) + a P_n(x,t) + (t + b) P_{n-1}(x,t)
\end{align*}
for $n \geq 2$. For $a = b = 0$, the polynomials are the Chebyshev polynomials of the second kind. For general $a \in \mf{R}$, $b \in \mf{R}_+$, such polynomials have been considered by a number of authors; see the discussion in pages 26--28 of \cite{Askey-Wilson}, and also \cite{Freeman}.

\begin{Defn}
Let $\mb{R}$ be a free cumulant generating function, and $\mb{U}(\mb{z})$ an $n$-tuple of non-com\-mu\-ta\-tive power series such that $U_i(\mb{z}) = z_i + \text{ higher-order terms}$. Multivariate free Sheffer polynomials are defined via their generating function
\[
H(\mb{x}, \mb{z})
= \bigl(1 - \mb{x} \cdot \mb{U}(\mb{z}) + \mb{R}(\mb{U}(\mb{z})) \bigr)^{-1}
= 1 + \sum_{\vec{u}} P_{\vec{u}}(\mb{x}) z_{\vec{u}}.
\]
\end{Defn}

\noindent
Define a linear operator $D_j$ (left partial derivative) on non-commutative power series via
\[
D_j w_{u(1), \ldots, u(n)} =
\begin{cases}
0, & u(1) \neq j, \\
w_{u(2), \ldots, u(n)}, & u(1) = j.
\end{cases}
\]

\begin{Thm}
Suppose that for a family of multivariate free Sheffer polynomials, $P_{\vec{u}}$ and $P_i$ are orthogonal for all $i$ and $\abs{\vec{u}} \geq 2$. In particular, this is the case if the polynomials are pseudo-orthogonal. Suppose also that the covariance matrix $R[X_j, X_i]$ is non-degenerate. Then
\[
\mb{U}(\mb{z}) = \mb{F}^{<-1>} \Bigl(\sum_{i=1}^n R[X_1, X_i] z_i, \ldots, \sum_{i=1}^n R[X_n, X_i] z_i \Bigr).
\]
Here $F_i(\mb{z}) = (D_i \mb{R}) (\mb{z}) - \state{X_i}$, and for an $n$-tuple of power series $\mb{F}$, $\mb{G} = \mb{F}^{<-1>}$ is the inverse of $\mb{F}$ under composition,
\[
F_i(\mb{G}(\mb{z})) = z_i.
\]

\end{Thm}

\noindent
In the classical case, the corresponding condition defines precisely the natural exponential families. The proof below is inspired by the one in \cite{Pommeret-Sheffer}.

\begin{proof}
Using the substitution $U_i(\mb{z}) = w_i (1 + \mb{M}(\mb{w}))$ leads to $\mb{R}(\mb{U}(\mb{z})) = \mb{M}(\mb{w})$. Then
\[
H(\mb{x}, \mb{z})
= \Bigl(1 - \mb{x} \cdot \mb{w} \bigl(1 + \mb{M}(\mb{w}) \bigr) + \mb{M}(\mb{w}) \Bigr)^{-1}
= (1 + \mb{M}(\mb{w}))^{-1} \bigl(1 - \mb{x} \cdot \mb{w} \bigr)^{-1}.
\]
It follows that $\state{H(\mb{X}, \mb{z})} = 1$. Suppose the conditions of the theorem are satisfied. By definition, $P_j = X_j - \state{X_j}$, so $P_j^\ast = P_j$. Then
\[
\state{P_j^\ast z_j H(\mb{X}, \mb{z})}
= \state{\sum_{\vec{u}} P_j P_{\vec{u}} z_j z_{\vec{u}}}
= \sum_{i=1}^n \state{P_j P_i} z_j z_i.
\]
The covariance of $X_j$ and $X_i$ is
\[
\state{(X_j - \state{X_j}) (X_i - \state{X_i})} = R[X_j, X_i].
\]
So
\[
\sum_{i=1}^n R[X_j, X_i] z_j z_i
= \state{P_j z_j H(\mb{X}, \mb{z})}
= \state{X_j z_j H(\mb{X}, \mb{z})} - \state{X_j} z_j
\]
and
\[
\state{X_j z_j H(\mb{X}, \mb{z})}
= \sum_{i=1}^n R[X_j, X_i] z_j z_i + \state{X_j} z_j.
\]
On the other hand, using the substitution above,
\begin{equation*}
\begin{split}
x_j z_j H(\mb{x}, \mb{z})
&= x_j z_j \bigl(1 - \mb{x} \cdot \mb{U}(\mb{z}) + \mb{R}(\mb{U}(\mb{z})) \bigr)^{-1} \\
&= x_j z_j \bigl(1 - \mb{x} \cdot \mb{w} (1 + \mb{M}(\mb{w})) + \mb{M}(\mb{w}) \bigr)^{-1} \\
&= z_j (1 + \mb{M}(\mb{w}))^{-1} x_j (1 - \mb{x} \cdot \mb{w})^{-1}
= z_j (1 + \mb{M}(\mb{w}))^{-1} \Bigl(x_j + \sum_{\vec{u}} x_j x_{\vec{u}} w_{\vec{u}} \Bigr),
\end{split}
\end{equation*}
and so
\[
\state{X_j z_j H(\mb{X}, \mb{z})}
= z_j (1 + \mb{M}(\mb{w}))^{-1} \Bigl(M[X_j] + \sum_{\vec{u}} M[X_j, X_{\vec{u}}] w_{\vec{u}} \Bigr).
\]
Since
\[
D_j \mb{M}(\mb{w})
= D_j \Bigl(\sum_{\vec{u}} M[X_{\vec{u}}] w_{\vec{u}} \Bigr)
= M[X_j] + \sum_{\vec{u}} M[X_j, X_{\vec{u}}] w_{\vec{u}},
\]
it follows that
\[
\state{X_j z_j H(\mb{X}, \mb{z})}
= z_j (1 + \mb{M}(\mb{w}))^{-1} D_j (\mb{M}(\mb{w})).
\]
Since $U_i(\mb{z}) = w_i (1 + \mb{M}(\mb{w}))$ and $\mb{R}(\mb{U}(\mb{z})) = \mb{M}(\mb{w})$,
\[
(1 + \mb{M}(\mb{w}))^{-1} D_{w_j}(\mb{M}(\mb{w})) = (D_{z_j} \mb{R})(\mb{U}(\mb{z})).
\]
We conclude that
\[
\state{X_j z_j H(\mb{X}, \mb{z})}
= z_j (D_j \mb{R})(\mb{U}(\mb{z}))
\]
and
\[
(D_j \mb{R})(\mb{U}(\mb{z})) = \sum_{i=1}^n R[X_j, X_i] z_i + R[X_j].
\]
Thus
\[
F_j(\mb{U}(\mb{z})) = \sum_{i=1}^n R[X_j, X_i] z_i.
\]
Since
\[
F_j(\mb{z})
= \sum_{i=1}^n R[X_j, X_i] z_i + \text{ higher-order terms}
\]
and the covariance matrix is non-degenerate, this series has an inverse under composition.
\end{proof}

\subsection{Free Kailath-Segall polynomials}
We only list three results in this section; all the other results from Section~\ref{Subsec:Kailath} have free analogs, which will be proven in greater generality in Section~\ref{Subsec:q-Kailath}. Free Kailath-Segall polynomials are the $q$-Kailath-Segall polynomials from that section for $q=0$.

\br
First, free Kailath-Segall polynomials have explicit expansions
\[
\KS{f_1, f_2,, \ldots, f_n}
= \sum_{\pi \in \Int(n)} \sum_{S \subset \Sing(\pi)} (-1)^{n - \abs{\pi} + \abs{S}} \prod_{B \in S} \Exp{f_B} \prod_{C \in S^c} X(f_C).
\]
Second, the product of free Kailath Segall polynomials can be expanded as
\[
\begin{split}
&\prod_{i=1}^k \KS{f_{u_i(1)}, f_{u_i(2)}, \ldots, f_{u_i(s(i))}} \\
&\quad = \sum_{\substack{\pi \in \NC(N) \\ \pi \wedge \pi_{s(1), s(2), \ldots, s(k)} = \hat{0}}} \sum_{\substack{S \subset \pi, \\ \Sing(S^c) = \emptyset}} \prod_{B \in S^c} \Exp{f_{(\vec{u}:B)}} W(f_{(\vec{u}:C)}: C \in S)
\end{split}
\]
(the proof is similar to the one in Proposition~\ref{Prop:Kailath-Segall}). Third, free Appell polynomials can be expressed in terms of the free Kailath-Segall polynomials.

\begin{Prop}
For $x_1 = x$,
\[
A^{(n)}(x)
= \sum_{k=1}^n \sum_{\vec{u} \in \Delta_{k, n}} W_{\vec{u}}(\mb{x}).
\]
\end{Prop}

\begin{proof}
The proof is similar to Proposition~\ref{Prop:Appell-KS}. We show that
\[
A^{(n)}(X) \cdot 1
= \sum_{k=1}^n \sum_{\vec{u} \in \Delta_{k, n}} \prod_{i=1}^k x_{u(i)}
\]
by induction, with the induction step being to verify that
\begin{equation*}
\begin{split}
A^{(n+1)}(X) \cdot 1
&= \sum_{k=1}^{n+1} \sum_{\vec{u} \in \Delta_{k, n+1}} \biggl[ x_1 \prod_{i=1}^k x_{u(i)} + x_{u(1) + 1} \prod_{i=2}^k x_{u(i)} + r_{u(1) + 1} \prod_{i=2}^k x_{u(i)} + r_1 \prod_{i=1}^k x_{u(i)} \biggr] \\
&\quad
- \sum_{j=0}^n r_{n+1-j} \sum_{k=1}^j \sum_{\vec{v} \in \Delta_{k, j}} \prod_{i=1}^k x_{v(i)}. \qedhere
\end{split}
\end{equation*}
\end{proof}

\section{$q$-interpolation}
\label{Section:Q}
\noindent
We saw in the preceding two sections that both the Appell and the Kailath-Segall polynomials are related to probability theory (they are martingale polynomials), and that they have analogs performing the same functions in free probability theory. On the other hand, there is also a very clear relationship between these classes and the third one of orthogonal polynomials. In the light of this, it is interesting to look at their $q$-deformations. For orthogonal polynomials, such deformations are given by the members of the Askey scheme of basic hypergeometric orthogonal polynomials \cite{Askey-Scheme}. On the other hand, a possible $q$-deformed probability theory was initiated by Bo\.{z}ejko and Speicher~\cite{BozSpeBM1}. We show that the Kailath-Segall polynomials can be extended to this context, with the same relation to the deformed probability theory and orthogonal polynomials as in the $q=0,1$ cases. As a consequence, we obtain combinatorial formulas for some basic hypergeometric orthogonal polynomials, and Wick product formulas for Fock representations. In contrast, the $q$-extension of the Appell polynomials is not satisfactory. We propose a definition in the one-variable case which fits well with the orthogonal polynomials picture (the analogs of the Meixner families are the Al-Salam and Chihara polynomials). However, this definition does not fit well with the $q$-deformed L\'{e}vy processes (see the Appendix), and an adequate definition in the multivariate case is missing.

\subsection{Notation}
Fix $q \in (-1, 1)$. A few standard pieces of $q$-notation are: $[0]_q = 0$,
\[
[n]_q = \sum_{i=1}^n q^i = \frac{1 - q^n}{1-q}
\]
for $n \geq 1$, $[0]_q! = 1$, $[n]_q! = \prod_{i=1}^n [i]_q$ for $n \geq 1$, $\binom{n}{k}_q = \frac{[n]_q!}{[k]_q! [n-k]_q!}$.

\br
Define the $q$-cumulants for a single measure via the relation
\[
m_n(\mu) = \sum_{\pi \in \Part(n)} q^{\rc{\pi}} \prod_{B \in \pi} r_{\abs{B}}(\mu).
\]
Define the cumulant generating function for this section to be
\[
\mb{R}_{\mu}(z) = \sum_{n=1}^\infty \frac{1}{[n]_q!} r_n(\mu) z^n.
\]
More generally, for a $n$-tuple of random variables, define the joint $q$-cumulants via
\begin{equation*}
M[X_{\vec{u}}]
= \sum_{\pi \in \NC(n)} q^{\rc{\pi}} \prod_{B \in \pi} R[X_{(\vec{u}:B)}].
\end{equation*}

\subsection{$q$-Appell polynomials}
There is a number of different possible definitions for $q$-Appell polynomials. For example, two $q$-deformations of the relation $\partial_x A_n = n A_{n-1}$, leading to such definitions, were considered in \cite{Al-Salam-q-Appell,Al-Salam-Hermite}. We prefer, instead, to use an interpolation between the recursion relations in the classical and the free cases.
\begin{Defn}
The $q$-Appell polynomials are defined via the recursion relation
\[
A_{n+1}(x) = x A_n(x) - \sum_{k=0}^n \binom{n}{k}_q r_{n+1-k} A_k(x),
\]
where $\set{r_k}$ are some $q$-cumulant sequence.
\end{Defn}

\begin{Prop}
\label{Prop:q-Appell-generating}
The generating function of $q$-Appell polynomials has the form
\[
\sum_{n=0}^\infty \frac{1}{[n]_q!} A_n(x) z^n
= H(x,z) = \prod_{k=0}^\infty \frac{1}{1 - (1-q) x z q^k + \mb{R}(z q^k) - \mb{R}(z q^{k+1})}.
\]
\end{Prop}

\begin{proof}
Denote $D_{q,z} z^n = [n]_q z^{n-1}$, so that
\[
D_{q,z}(f) = \frac{f(z) - f(q z)}{z - q z}
\]
is the standard $q$-derivative operator. By definition,
\[
\frac{1}{[n]_q!} A_{n+1} z^{n+1}
= x \frac{1}{[n]_q!} A_n z^{n+1} - \sum_{k=0}^n \frac{1}{[n-k]_q!} r_{n+1-k} \frac{1}{[k]_q!} A_k z^{n+1}.
\]
Divide by $[n+1]_q$ to get
\[
\frac{1}{[n+1]_q!} A_{n+1} z^{n+1}
= x \frac{1}{[n+1]_q!} A_n z^{n+1} - \frac{1}{[n+1]_q} \sum_{k=0}^n \frac{1}{[n-k]_q!} r_{n+1-k} \frac{1}{[k]_q!} A_k z^{n+1}.
\]
Apply $D_{q,z}$ to get
\[
D_{q,z} \frac{1}{[n+1]_q!} A_{n+1} z^{n+1}
= x \frac{1}{[n]_q!} A_n z^n - \sum_{k=0}^n \frac{1}{[n-k]_q!} r_{n+1-k} z^{n-k} \frac{1}{[k]_q!} A_k z^{k}.
\]
So
\[
D_{q,z} H = x H - D_{q,z}(\mb{R}) H.
\]
It follows that
\[
H(x,z)
= \frac{1}{1 - (1-q) x z + \mb{R}(z) - \mb{R}(q z)} H(x, q z). \qedhere
\]
\end{proof}

\begin{Ex}
\label{Ex:q-Hermite}
Since orthogonal polynomials satisfy three-term recursion relations, the only orthogonal polynomials among the $q$-Appell are those with $r_k = 0$ for $k > 2$. Thus $\mb{R}(z) = a z + b z^2$. By adding a constant to $x$ and re-scaling the polynomials, we may take $r_k = \delta_{2 k} t $, $\mb{R}(z) = \frac{t}{1+q} z^2$, in which case
\[
H(x,z) = \frac{1}{1 - (1-q) x z + (1 - q) t z^2} H(x, q z).
\]
Up to further re-scaling, this is the generating function for the continuous (Rogers) $q$-Hermite polynomials.
\end{Ex}

\begin{Ex}
Consider the polynomial family defined by the recursion
\[
x P_n = P_{n+1} + (t - [n]_q) P_n + [n]_q x P_{n-1}.
\]
It has the generating function satisfying
\[
H(x,t,z) = \frac{1}{1 + (1-q) \frac{z}{1-z} t - (1-q) z x} H(x,t,q z).
\]
So this is a family of $q$-Appell polynomials for $\mb{R}(z) = t \sum_{n=0}^\infty \frac{1}{[n]_q} z^n$, $r_n = t [n-1]_q!$. As $q \rightarrow 0$, $H$ converges to
\[
\frac{1}{1 + \frac{z}{1-z} t - z x},
\]
the generating function for the free Appell polynomials for the free Poisson process. On the other hand, for $q \rightarrow 1$, $H$ converges to
\[
e^{x z} (1-z)^t,
\]
the generating function for the Appell polynomials for the Gamma process. This can be considered as an extension of a well-known property that for $t=1$, a square of a normal random variable has a Gamma distribution, while a square of a semicircular random variable has a free Poisson distribution.
\end{Ex}

\begin{Prop}
Let $\set{A_n}$ be the $q$-Appell polynomials for the $q$-cumulant generating function $\mb{R}$. Then
\[
\partial_x A_n(x) = \sum_{k=0}^{n-1} \binom{n}{k}_q A_k(x) A_{n-k-1}^0(x),
\]
where $\set{A_n^0}$ are defined via their generating function
\[
\sum_{n=0}^\infty \frac{1}{[n]_q!} A_n^0(x) z^n
= \frac{1}{1 - (1-q) x z + \mb{R}(z) - \mb{R}(q z)}.
\]
\end{Prop}

\begin{proof}
By Proposition~\ref{Prop:q-Appell-generating},
\begin{equation*}
\begin{split}
\partial_x H(x,z)
&= \partial_x \prod_{k=0}^\infty \frac{1}{1 - (1-q) x z q^k + \mb{R}(z q^k) - \mb{R}(z q^{k+1})} \\
&= \sum_{k=0}^\infty H(x,z) \frac{1}{1 - (1-q) x z q^k + \mb{R}(z q^k) - \mb{R}(z q^{k+1})} (1-q) z q^k \\
&= \sum_{n=0}^\infty \frac{1}{[n]_q!} A_n(x) z^n \sum_{k=0}^\infty \sum_{l=0}^\infty \frac{1}{[l]_q!} A_l^0(x) (z q^k)^l (1-q) z q^k \\
&= \sum_{n=0}^\infty \frac{1}{[n]_q!} A_n(x) z^n \sum_{l=0}^\infty \frac{1}{[l]_q!} A_l^0(x) \frac{1-q}{1 - q^{l+1}} z^{l+1} \\
&= \sum_{n=0}^\infty \sum_{l=0}^\infty \frac{1}{[n]_q! [l+1]_q!} A_n(x) A_l^0(x) z^{n + l + 1}. \qedhere
\end{split}
\end{equation*}
\end{proof}

\begin{Ex}
If $\set{P_n(x)}$ are the continuous $q$-Hermite polynomials, then
\[
\partial_x P_n(x) = \sum_{k=0}^{n-1} \frac{[n]_q!}{[n-k]_q! [k]} P_k(x,t) Q_{n-k-1}(x),
\]
where
\[
\sum_{n=0}^\infty Q_n(x)
= \frac{1}{1 - (1-q) x z + (1-q) z^2}.
\]
So $Q_n(x) z^n = U_n(\sqrt{1-q} x) (\sqrt{1-q} z)^n$, where $U_n$ are the Chebyshev polynomials of the second kind. Therefore
\[
\partial_x P_n(x) = \sum_{k=0}^{n-1} \frac{[n]_q!}{[n-k]_q! [k]} (\sqrt{1-q})^{k-1} P_k(x) U_{n-k-1}(\sqrt{1-q} \,x).
\]
\end{Ex}

\subsection{Al-Salam-Chihara polynomials}
In this section we consider a $q$-deformation of the Meixner families. One such deformation has been considered in \cite{Ismail-Operator-Calculus}. It is based on the exponential function for the operator calculus for the Askey-Wilson operator \cite{Ismail-Zhang}. Previously, Al-Salam showed that under this approach, the unique orthogonal Appell family also consists of (multiples of) Rogers $q$-Hermite polynomials \cite{Al-Salam-Hermite}. Under the more elementary approach of \cite{Al-Salam-q-Appell}, the unique orthogonal Appell family are the Al-Salam-Carlitz polynomials.

\br
However, we have an extra requirement to put on our families. In addition to the correct limiting behavior as $q \rightarrow 1$, we also require a correct limiting behavior as $q \rightarrow 0$. Specifically, in this limit we should obtain the free Meixner families.

\br
Define orthogonal polynomials by the recursion relation
\begin{equation}
\label{q-recursion}
x P_n = P_{n+1} + a [n]_q P_n + [n]_q (t + b [n-1]_q) P_{n-1}
\end{equation}
for $a \in \mf{R}$, $b \in \mf{R}_+$ (cf.\ \cite[Remark 6]{AnsMeixner}). Then the free Meixner families are these families of polynomials for $q=0$, while the classical Meixner families are these families for $q=1$.

\begin{Lemma}
Up to re-scaling, the polynomials defined by recursion~\eqref{q-recursion} are the Al-Salam-Chihara polynomials.
\end{Lemma}

\begin{proof}
Define the generating function $H(x,t,z) = \sum_{n=0}^\infty \frac{1}{[n]_q!} P_n(x,t) z^n$. Then
\begin{equation*}
\begin{split}
x z (H(x,t,z) - 1)
&= \frac{1}{1-q} [H - P_1(x,t) z - H(x,t,q z) + P_1(x,t) q z] \\
&\quad+ \frac{a z}{1-q} [H(x,t,z) - H(x,t,q z)] \\
&\quad+ t z^2 H(x,t,z) + \frac{b z^2}{1-q} [H(x,t,z) - H(x,t,q z)].
\end{split}
\end{equation*}
Since $P_1(x,t) = x$,
\[
H(x,t,z) \left( \frac{1}{1-q} - x z + \frac{a z}{1-q} + t z^2 + \frac{b z^2}{1-q} \right)
= H(x,t,q z) \left( \frac{1}{1-q} + \frac{a z}{1-q} + \frac{b z^2}{1-q} \right).
\]
Therefore
\begin{equation*}
\begin{split}
H(x,t,z)
&= \frac{(1 + a z + b z^2)}{(1 + a z + b z^2) + (1-q) (t z^2 - x z)} H(x,t,q z) \\
&= \frac{(1 + a z + b z^2)}{1 + [a - (1-q) x]z + [b + (1-q) t]z^2} H(x,t,q z) \\
&= \frac{1}{1 + \frac{1-q}{(1 + a z + b z^2)} (t z^2 - x z)} H(x,t,q z).
\end{split}
\end{equation*}
Affine transformations $w = \sqrt{b + (1-q)t} z$, $y = \frac{(1-q) x - a}{2 \sqrt{b + (1-q)t}}$ bring $H$ to the standard form of the generating function for the Al-Salam-Chihara polynomials \cite[3.8]{Askey-Scheme}
\[
H(y,w;a',b') = \frac{(1-a' w)(1-b' w)}{1 - 2 y w + w^2} H(y,q w;a',b'),
\]
where $a', b'$ are the roots of the quadratic polynomial $z^2 + \frac{a}{\sqrt{b + (1-q)t}} z + b$. 
\end{proof}

\noindent
In particular, the Hermite case $a = b = 0$ corresponds to the Rogers (continuous) $q$-Hermite polynomials, while the Charlier case $b = 0$ corresponds to what are usually called the continuous big $q$-Hermite polynomials. Note also that the $q$-Krawtchouk polynomials defined in \cite{SaiKraw} are, up to a shift, of this form, with $a = 1-2p$, $b = -p(1-p)$, $2t = N p(1-p)$, except that in this case, as expected, $b$ is negative.

\br
Al-Salam-Chihara polynomials were defined in \cite{Al-Salam-Chihara} as all polynomials other than the Meixner families characterized by a certain convolution property. Their measure of orthogonality was found explicitly in \cite{Askey-Ismail}. The interpretation of these polynomials as $q$-analogs of the Meixner families was explicitly conjectured in \cite{Andrews-Askey}, and proved in \cite{Al-Salam-Pollaczek}. We found the following proof independently, and in our particular case it is also somewhat simpler. For a more interesting characterization using stochastic processes, see \cite{Bryc-Meixner}.

\begin{Thm}
\label{Thm:Characterization}
Suppose the generating function $H(x,z) = \sum_{n=0}^\infty \frac{1}{[n]_q!} P_n(x) z^n$ of the monic orthogonal polynomials defined by the recursion relation
\[
x P_n(x) = P_{n+1}(x) + \alpha_n P_n(x) + \beta_n P_{n-1}(x)
\]
has the form
\[
H(x,z) = F(z) \prod_{k=0}^\infty \frac{1}{1 - (1-q) U(q^k z) x} = \frac{F(z)}{F(q z)} \frac{1}{1 - (1-q) U(z) x} H(x, q z),
\]
where $F, U$ are formal power series with $F(z) = 1 + \text{ higher-order}$ terms, $U(z) = z + \text{ higher-order}$ terms. Then the polynomials are a re-scaling of the Al-Salam-Chihara polynomials, with
\[
\alpha_n = \alpha_0 + c_1 [n]_q
\]
and
\[
\beta_n = [n]_q (\beta_1 + c_2 [n-1]_q).
\]
\end{Thm}

\begin{proof}
By assumption,
\[
(1 - (1-q) U(z) x) H(x,z) = \frac{F(z)}{F(q z)} H(x,q z).
\]
Therefore
\begin{equation}
\label{x-H}
x H(x,z) = \frac{1}{(1-q) U(z)} \left( H(x,z) - \frac{F(z)}{F(q z)} H(x,q z) \right).
\end{equation}
Define the lowering operator $D$ on $\mf{R}[x]$ by $D(P_n) = [n]_q P_{n-1}$ and extend linearly. Then $D(H) = z H$ and equation~\eqref{x-H} implies
\begin{equation}
\label{A-H}
\begin{split}
D(x H)(x,z)
&= \frac{1}{(1-q) U(z)} \left( z H(x,z) -  \frac{F(z)}{F(q z)} q z H(x,q z) \right) \\
&= z x H(x,z) + \frac{z}{U(z)} \frac{F(z)}{F(q z)} H(x,q z).
\end{split}
\end{equation}
Expand $\frac{z}{U(z)} \frac{F(z)}{F(q z)}$ into the formal power series $\sum_{n=0}^\infty c_n z^n$, with $c_0 = 1$. Equation~\eqref{A-H} says
\begin{equation}
\label{A-P_n}
D(x P_n)(x)
= [n]_q x P_{n-1}(x) + \sum_{k=0}^n c_k \frac{[n]_q!}{[n-k]!} q^{n-k} P_{n-k}(x).
\end{equation}
Applying the lowering operator to the recursion relation gives
\begin{equation}
\label{A-P_n-2}
D(x P_n)(x)
= [n+1]_q P_n(x) + \alpha_n [n]_q P_{n-1}(x) + \beta_n [n-1]_q P_{n-2}(x).
\end{equation}
The recursion relation for $n-1$ is
\[
x P_{n-1}(x) = P_n(x) + \alpha_{n-1} P_{n-1}(x) + \beta_{n-1} P_{n-2}(x).
\]
Multiplying it by $[n]_q$ gives
\begin{equation}
\label{Recursion-(n-1)}
[n]_q x P_{n-1}(x) = [n]_q P_n(x) + [n]_q \alpha_{n-1} P_{n-1}(x) + [n]_q \beta_{n-1} P_{n-2}(x).
\end{equation}
Combining Equations \eqref{A-P_n}, \eqref{A-P_n-2}, and \eqref{Recursion-(n-1)}, we obtain
\begin{multline*}
[n]_q P_n + \alpha_{n-1} [n]_q P_{n-1} + \beta_{n-1} [n]_q P_{n-2} + \sum_{k=0}^n c_k \frac{[n]_q!}{[n-k]_q!} q^{n-k} P_{n-k} \\
- [n+1]_q P_n - \alpha_n [n]_q P_{n-1} - \beta_n [n-1]_q P_{n-2} = 0.
\end{multline*}
Collecting coefficients of $P_{n-1}$ gives
\begin{equation}
\label{alpha}
\alpha_{n-1} [n]_q + c_1 [n]_q q^{n-1} - \alpha_n [n]_q = 0;
\end{equation}
collecting coefficients of $P_{n-2}$ gives
\begin{equation}
\label{beta}
\beta_{n-1} [n]_q + c_2 [n]_q [n-1]_q q^{n-2} - \beta_n [n-1]_q = 0,
\end{equation}
and collecting coefficients of $P_{n-k}$ for $k>2$ gives
\[
c_k = 0 \text{ for } k > 2.
\]
Equation~\eqref{alpha} gives
\[
\alpha_n - \alpha_{n-1} = c_1 q^{n-1},
\]
while equation~\eqref{beta} gives
\[
\frac{\beta_n}{[n]_q} - \frac{\beta_{n-1}}{[n-1]_q} = c_2 q^{n-2}.
\]
We conclude that $\alpha_n = \alpha_0 + c_1 [n]_q$ and $\beta_n = [n]_q (\beta_1 + c_2 [n-1]_q)$.
\end{proof}

\br
Thus it appears reasonable to define the $q$-Sheffer polynomials via their generating function
\begin{equation}
\label{q-Sheffer}
\sum_{n=0}^\infty \frac{1}{[n]_q!} P_n(x,t) z^n
= \prod_{k=0}^\infty \frac{1}{1 - (1-q) x U(z q^k) + \mb{R}(U(z q^k)) - \mb{R}(U(z q^{k+1}))}.
\end{equation}

\subsection{$q$-Kailath-Segall polynomials}
\label{Subsec:q-Kailath}
The origin of these polynomials is in the $q$-L\'{e}vy processes defined in \cite{AnsQCum,AnsQLevy}.

\begin{Defn}
Let $\mc{A}_0$ be a complex star-algebra without identity, and $\Exp{\cdot}$ a star-linear functional on it. Let $\mc{A}$ be the complex unital star-algebra generated by non-commuting symbols $\set{X(f): f \in \mc{A}_0^{sa}}$ (and $1$) subject to the linearity relations
\[
X(\alpha f + \beta g) = \alpha X(f) + \beta X(g).
\]
Equivalently, $\mc{A}$ is the tensor algebra of $\mc{A}_0$. The star-operation on it is determined by the requirement that all $X(f), f \in \mc{A}_0^{sa}$ are self-adjoint. For such $f_i$, define the $q$-Kailath-Segall polynomials by $\KS{f} = X(f) - \Exp{f}$ and
\begin{equation}
\label{q-KS}
\begin{split}
\KS{f, f_1, f_2, \ldots, f_n}
& = X(f) \KS{f_1, f_2, \ldots, f_n}
- \sum_{i=1}^n q^{i-1} \Exp{f f_i} \KS{f_1, \ldots, \hat{f}_i, \ldots, f_n} \\
&\quad - \sum_{i=1}^n q^{i-1} \KS{f f_i, \ldots, \hat{f}_i, \ldots, f_n}
- \Exp{f} \KS{f_1, f_2, \ldots, f_n}.
\end{split}
\end{equation}
This map has a $\mf{C}$-linear extension, so that each $W$ is really a multi-linear map from $\mc{A}_0$ to $\mc{A}$.
\end{Defn}

\noindent
In the particular case $\mc{A}_0 = \mf{C}_0[x]$ (polynomials without constant term), we may denote $x_i = X(x^i)$. The functional can be taken to be the $q$-cumulant functional of a measure $\mu$, $\Exp{x^i} = r_i(\mu)$. Then
\[
W_{\vec{u}}(\mb{x}) = \KS{x^{u(1)}, x^{u(2)}, \ldots, x^{u(n)}}
\]
are multivariate polynomials in $\set{x_i: i \in \mf{N}}$.

\begin{Prop}
The Kailath-Segall formula~\eqref{KS-formula} takes the following form.
\begin{equation*}
W(\underbrace{f, f, \ldots, f}_{n+1}) 
= \sum_{k=0}^n (-1)^k \frac{[n]_q!}{[n-k]_q!} X(f^{k+1}) W(\underbrace{f, f, \ldots, f}_{n-k}) 
- \Exp{f} W(\underbrace{f, f, \ldots, f}_n).
\end{equation*}
\end{Prop}

\noindent
The proposition can be proven directly by induction, or deduced from Theorem~\ref{Thm:Q-expansions}(b).

\br
Many of the following formulas appear in \cite{Effros-Popa} in the $q$-Gaussian case $\mc{A}_0 = \mf{C}_0 \langle \mb{x} \rangle$ and $\Exp{x_{\vec{u}}} = \delta_{\abs{\vec{u}}, 2}$.

\begin{Thm}
\label{Thm:Q-expansions}
The following expansions hold.
\begin{enumerate}
\item
A monomial in $X(f)$'s can be expanded in terms of the $q$-Kailath-Segall polynomials:
\begin{equation*}
X(f_1) X(f_2) \ldots X(f_n) =
\sum_{\pi \in \Part(n)} \sum_{S \subset \pi} q^{\rc{S, \pi}} \prod_{B \in S^c} \Exp{f_B} \KS{f_C: C \in S}.
\end{equation*}
\item
For a permutation $\sigma \in \Sym(n)$, write its standard cycle decomposition as
\[
\sigma = \bigl(u(1,1), \ldots, u(1, s_1) \bigr)  \ldots \bigl(u(k,1), \ldots, u(k, s_k) \bigr).
\]
Here $u(i,1) = \min_j \set{u(i,j)}$, and the cycles of $\sigma$ are ordered according to the order of their smallest elements, $u(1,1) < u(2,1) < \ldots < u(k, 1)$. Let $\Sing(\sigma)$ be the one-element cycles of $\sigma$, and let $s(\sigma)$  the number of inversions of the permutation
\[
F(\sigma) =
\left(
\begin{matrix}
1 & \ldots & n \\
u(1,1) & \ldots & u(k, s_k)
\end{matrix}
\right)
\]
($F$ is almost, but not quite, the fundamental transformation of Foata \cite{Foata}). Finally, for a subset $S \subset \Sing(\sigma)$, denote
\[
X_{(S,\sigma)}(f_1, f_2, \ldots, f_n) = \prod_{i: (u(i,1), \ldots, u(i, s_i)) \in S^c} X(f_{u(i,1)} f_{u(i,2)} \ldots f_{u(i, s_i)}).
\]
Note that in fact, $X_{(S, \sigma)}$ depends only on $\set{f_i: \set{i} \not \in S}$. With this notation, the $q$-Kailath-Segall polynomials are
\[
\KS{f_1, f_2, \ldots, f_n}
= \sum_{\sigma \in \Sym(n)} \sum_{S \subset \Sing(\sigma)} (-1)^{n - \cyc(\sigma) + \abs{S}} q^{s(\sigma)} \prod_{i: \set{i} \in S} \Exp{f_i}  X_{(S,\sigma)}(f_1, f_2, \ldots, f_n).
\]
\end{enumerate}
\end{Thm}

\begin{proof}
Part (a) was proven in the appendix of \cite{AnsQLevy} using the Fock space representation; one can also use the defining recursion relation. For part (b), by definition,
\begin{equation*}
\begin{split}
\KS{f, f_1, f_2, \ldots, f_n}
& = X(f) \KS{f_1, f_2, \ldots, f_n}
- \sum_{i=1}^n q^{i-1} \Exp{f f_i} \KS{f_1, \ldots, \hat{f}_i, \ldots, f_n} \\
&\quad - \sum_{i=1}^n q^{i-1} \KS{f f_i, \ldots, \hat{f}_i, \ldots, f_n}
- \Exp{f} \KS{f_1, f_2, \ldots, f_n}.
\end{split}
\end{equation*}
Using the defining recursion and induction on $n$, this is
\begin{align*}
& \qquad X(f) \sum_{\tau \in \Sym(n)} \sum_{S \subset \Sing(\tau)} (-1)^{n - \cyc(\tau) + \abs{S}} q^{s(\tau)} \prod_{j \in S} \Exp{f_j} X_{(S,\tau)}(f_1, \ldots, f_n) \\
& - \sum_{i=1}^n q^{i-1} \Exp{f f_i} \sum_{\tau \in \Sym\left(\set{1, \hat{i}, n}\right)} \sum_{S \subset \Sing(\tau)} (-1)^{n - 1 - \cyc(\tau) + \abs{S}} q^{s(\tau)} \prod_{j \in S} \Exp{f_j} X_{(S,\tau)}(f_1, \ldots, \hat{f}_i, \ldots, f_n) \\
& - \sum_{i=1}^n q^{i-1} \sum_{\tau \in \Sym \left( \set{0, 1, \hat{i}, n} \right)} \sum_{\substack{S \subset \Sing(\tau) \\ 0 \not \in S}} (-1)^{n - \cyc(\tau) + \abs{S}} q^{s(\tau)} \prod_{j \in S} \Exp{f_j} X_{(S,\tau)}(f f_i, f_1, \ldots, \hat{f}_i, \ldots, f_n) \\
& - \sum_{i=1}^n q^{i-1} \sum_{\tau \in \Sym \left( \set{0, 1, \hat{i}, n} \right)} \sum_{\substack{S \subset \Sing(\tau) \\ 0 \in S}} (-1)^{n - \cyc(\tau) + \abs{S}} q^{s(\tau)} \Exp{f f_i} \prod_{\substack{j \in S \\ j \neq 0}} \Exp{f_j} \\
&\qquad \times X_{(S,\tau)}(f f_i, f_1, \ldots, \hat{f}_i, \ldots, f_n) \\
& - \Exp{f} \sum_{\tau \in \Sym(n)} \sum_{S \subset \Sing(\tau)} (-1)^{n - \cyc(\tau) + \abs{S}} q^{s(\tau)} \prod_{j \in S} \Exp{f_j} X_{(S,\tau)}(f_1, \ldots, f_n)
\end{align*}
Since in the fourth term, $0 \in \Sing(\tau)$, the second and the fourth terms cancel. So we obtain
\begin{equation*}
\begin{split}
& \quad \sum_{\tau \in \Sym(n)} \sum_{S \subset \Sing(\tau)} (-1)^{(n+1) - (\cyc(\tau)+1) + \abs{S}} q^{s(\tau)} \prod_{j \in S} \Exp{f_j} X(f) X_{(S,\tau)}(f_1, \ldots, f_n) \\
& + \sum_{i=1}^n \sum_{\tau \in \Sym \left( \set{0, 1, \hat{i}, n} \right)} \sum_{\substack{S \subset \Sing(\tau) \\ 0 \not \in S}} (-1)^{(n+1) - \cyc(\tau) + \abs{S}} q^{s(\tau)+i-1} \prod_{j \in S} \Exp{f_j}  X_{(S,\tau)}(f f_i, f_1, \ldots, \hat{f}_i, \ldots, f_n) \\
& + \sum_{\tau \in \Sym(n)} \sum_{S \subset \Sing(\tau)} (-1)^{(n+1) - (\cyc(\tau)+1) + (\abs{S}+1)} q^{s(\tau)} \Exp{f} \prod_{j \in S} \Exp{f_j} X_{(S,\tau)}(f_1, \ldots, f_n).
\end{split}
\end{equation*}
The three terms in the preceding equation correspond to pairs $(S, \sigma)$, $\sigma \in \Sym \left(\set{0, 1, \ldots, n} \right)$ such that
\begin{align*}
& 0 \in \Sing(\sigma), 0 \not \in S & \text{(first term)}, \\
& 0 \not \in \Sing(\sigma), \sigma(0) = i & \text{(second term)}, \\
& 0 \in S \subset \Sing(\sigma) & \text{(third term)}.
\end{align*}
It remains to match up the powers of $q$. Suppose that
\[
\sigma = (u(1,1), \ldots, u(1, s_1))  \ldots (u(k,1), \ldots, u(k, s_k)),
\]
$s(\sigma) = i(F(\sigma))$ for
\[
F(\sigma) =
\left(
\begin{matrix}
0 & \ldots & n \\
u(1,1) & \ldots & u(k, s_k)
\end{matrix}
\right).
\]
In the first and the third terms,
\[
F(\sigma) =
\left(
\begin{matrix}
0 & 1 & \ldots & n \\
0 & u(2,1) & \ldots & u(k, s_k)
\end{matrix}
\right),
\]
and $i(F(\sigma)) = i(F(\tau))$ for
\[
F(\tau) = \left(
\begin{matrix}
1 & \ldots & n \\
u(2,1) & \ldots & u(k, s_k)
\end{matrix}
\right).
\]
In the second term,
\[
F(\sigma) =
\left(
\begin{matrix}
0 & 1 & \ldots & n \\
0 & i & \ldots & u(k, s_k)
\end{matrix}
\right),
\]
and $i(F(\sigma)) = i(F(\tau)) + i - 1$ for
\[
F(\tau) = \left(
\begin{matrix}
0 & 1 & \ldots & n \\
0 & u(1,3) \text{ or } u(2,1) & \ldots & u(k, s_k)
\end{matrix}
\right).
\]
\end{proof}

\noindent
Define a unital linear functional $\phi$ on $\mc{A}$ by
\[
\state{\KS{f_1, f_2, \ldots, f_n}} = 0
\]
for all $n>0$, $f_1, \ldots, f_n \in \mc{A}$.

\begin{Cor}
\label{Cor:q-functional}
The functional $\phi$ is given by
\begin{equation}
\label{Q-correlation}
\state{X(f_1) X(f_2) \ldots X(f_n)} =
\sum_{\pi \in \Part(n)} q^{\rc{\pi}} \prod_{B \in \pi} \Exp{f_B}.
\end{equation}
That is, with respect to $\phi$, the $q$-cumulants of such an $n$-tuple are
\[
R[X(f_1), X(f_2), \ldots, X(f_n)] = \Exp{f_1 f_2 \ldots f_n}.
\]
\end{Cor}

\begin{Cor}
\label{Cor:Q-linearizations}
The linearization coefficient for the $q$-Kailath-Segall polynomials is
\[
\state{\prod_{i=1}^k \KS{f_{u_i(1)}, f_{u_i(2)}, \ldots, f_{u_i(s(i))}}}
= \sum_{\substack{\pi \in \Part(N) \\ \pi \wedge \pi_{s(1), s(2), \ldots, s(k)} = \hat{0}, \\ \Sing(\pi) = \emptyset}} q^{\rc{\pi}} \prod_{B \in \pi} \Exp{f_{(\vec{u}:B)}}.
\]
In other words, it is the sum of partitioned $q$-cumulants over all inhomogeneous partitions with no singletons, with weights $q^{\rc{\pi}}$.
\end{Cor}

\subsection{Fock space realization}
\label{Section:KS-Fock-space}
The $q$-Kailath-Segall polynomials form a monic polynomial family in infinitely many variables. So we can construct a Fock space for them as in Section~\ref{Section:Fock}; this space will be infinite-dimensional. Instead, we will construct the Fock space directly from the multi-linear maps $W$. As a vector space, it will be the space of all polynomials in elements of $\mc{A}_0$, modulo the linearity relations. Equivalently, it is $\bigoplus_{n=0}^\infty \mc{A}_0^{\otimes n}$. The induced inner product is determined by
\begin{equation}
\label{q-inner-product}
\ip{\bigotimes_{i=1}^k f_{u(i)}}{\bigotimes_{j=1}^n f_{v(j)}} = \delta_{n k} \sum_{\sigma \in \Sym(n)} q^{i(\sigma)} \prod_{j=1}^n \Exp{f_{u(j)} f_{v(\sigma(j))}}.
\end{equation}
If $\Exp{\cdot}$ is a faithful state on $\mc{A}_0$, for $q \in (-1, 1)$ this inner product is positive definite \cite{BozSpeBM1}. Hence $\phi$ is a state.

\br
For the classical case $q=1$, this inner product is only positive semi-definite even if $\Exp{\cdot}$ is positive definite. The quotient by the kernel of $\phi$ gives precisely the symmetric Fock space.

\br
From the defining recursion relation, the action of the operator $X(f)$ is
\begin{equation*}
\begin{split}
X(f) (f_1 \otimes f_2 \otimes \ldots \otimes f_n)
&= f \otimes f_1 \otimes f_2 \otimes \ldots \otimes f_n \\
&\quad+ \sum_{i=1}^n q^{i-1} \Exp{f f_i} f_1 \otimes \ldots \otimes \hat{f}_i \otimes \ldots \otimes f_n \\
&\quad+ \sum_{i=1}^n q^{i-1} f f_i \otimes \ldots \otimes \hat{f}_i \otimes \ldots \otimes f_n \\
&\quad+ \Exp{f} f_1 \otimes f_2 \otimes \ldots \otimes f_n.
\end{split}
\end{equation*}
Thus it is a sum of a creation operator
\begin{align*}
a^\ast(f): f_1 \otimes f_2 \otimes \ldots \otimes f_n &\mapsto f \otimes f_1 \otimes f_2 \otimes \ldots \otimes f_n, \\
\intertext{an annihilation operator}
a(f): f_1 \otimes f_2 \otimes \ldots \otimes f_n &\mapsto \sum_{i=1}^n q^{i-1} \Exp{f f_i} f_1 \otimes \ldots \otimes \hat{f}_i \otimes \ldots \otimes f_n, \\
\intertext{a preservation operator}
p(f): f_1 \otimes f_2 \otimes \ldots \otimes f_n &\mapsto \sum_{i=1}^n q^{i-1} f f_i \otimes \ldots \otimes \hat{f}_i \otimes \ldots \otimes f_n, \\
\intertext{and a scalar operator}
f_1 \otimes f_2 \otimes \ldots \otimes f_n &\mapsto \Exp{f} f_1 \otimes f_2 \otimes \ldots \otimes f_n.
\end{align*}

\br
Now consider the case when $\Exp{\cdot}$ is positive semi-definite but not faithful. For simplicity, we will also assume that it has the trace property, $\Exp{fg} = \Exp{gf}$ for all $f, g \in \mc{A}_0$. Note that this does not imply that $\phi$ is a trace. Denote by
\[
I_{\Exp{\cdot}} = \set{f \in \mc{A}_0: \Exp{f^\ast f} = 0}
\]
the kernel of $\Exp{\cdot}$. Denote by $\mc{H}_0$ the Hilbert space obtained by completing the vector space $\mc{A}/I_{\Exp{\cdot}}$ with respect to the norm induced by $\norm{f} = \Exp{f^\ast f}$, with the induced inner product. Suppose $I_{\Exp{\cdot}}$ is in fact an ideal, so that
\begin{equation}
\label{Cstar}
\forall f \in I_{\Exp{\cdot}} \forall g \in \mc{A}_0, \Exp{f^\ast g^\ast g f} = 0.
\end{equation}
This is the case, for example, when $\mc{A}_0$ is a $C^\ast$-algebra. Note that
\[
\abs{\Exp{f g}} \leq \sqrt{\Exp{f^\ast f} \Exp{g^\ast g}},
\]
and so $\Exp{\mc{A}_0 I_{\Exp{\cdot}}} = 0$. It follows that the Fock representation of $\mc{A}$ factors through to the representation on $\bigoplus_{k=0}^\infty \mc{H}_0^{\otimes k}$ with the inner product induced from \eqref{q-inner-product}.

\br
Note that $\mc{A}_0 / I_{\Exp{\cdot}}$ is an algebra. In the examples below, we will observe the following situations. Let $\mf{C} 1 \oplus \mc{A}_0$ be the standard unital extension of the non-unital algebra $\mc{A}_0$.
\begin{enumerate}
\item
If the functional $\Exp{\cdot}$ can be extended to $\mf{C} 1 \oplus \mc{A}_0$ in a positive way, then $\abs{\Exp{f}}^2 \leq \Exp{f^\ast f} \Exp{1}$, and so $\Exp{I_{\Exp{\cdot}}} = 0$. Thus the action of $X(f)$ in the Fock representation depends only on the class of $f$ in $\mc{A}_0 / I_{\Exp{\cdot}}$. So we may replace $\mc{A}_0$ by $\mc{A}_0 / I_{\Exp{\cdot}}$ throughout.
\item
If the functional $\Exp{\cdot}$ has no such positive extension, it is natural to take $X(f)$ to be the sum of only three operators, $X(f) = a^\ast(f) + a(f) + p(f)$. The formulas for the $q$-Kailath-Segall polynomials need to be modified accordingly; this approach was used in \cite{AnsQLevy}. Under this construction, again the action of $X(f)$ in the Fock representation depends only on the class of $f$ in $\mc{A}_0 / I_{\Exp{\cdot}}$, and so we may replace $\mc{A}_0$ by $\mc{A}_0 / I_{\Exp{\cdot}}$ throughout.
\item
A particular case of part (b) is when
\begin{equation}
\label{Gaussian}
\mc{A}_0 \mc{A}_0 \subset  I_{\Exp{\cdot}}.
\end{equation}
In this case in the representation of $\mc{A}_0 / I_{\Exp{\cdot}}$, $X(f) = a^\ast(f) + a(f)$.
\end{enumerate}

\begin{Ex}
Throughout this example, $\mc{A}_0 = \mf{C}_0[x]$.
\begin{enumerate}
\item
Let $\nu$ be a positive measure on $\mf{R}$ all of whose moments are finite, and $\Exp{x^i} = m_i(\nu)$ for $i \geq 1$. Then we are in the context of part (a) of the preceding alternative, and so may replace $\mc{A}_0$ by $\mc{A}_0 / I_{\Exp{\cdot}}$ throughout. This is the case of compound $q$-Poisson distributions.

\br
More specifically, let $\Exp{x^i} = t$ for $i \geq 1$, so that $\nu = t \delta_1$. Then
\[
I_{\Exp{\cdot}} = \set{P: P(1) = 0}
\]
is an ideal. This is the $q$-Poisson case. In $\mc{A}_0 / I_{\Exp{\cdot}}$, $x^i = x$ for all $i$. So the $q$-Kailath-Segall polynomials in this case can be considered as single-variable polynomials. Specifically, they are the continuous big $q$-Hermite polynomials; the recursion relation for their centered version is a particular case of equation~\eqref{q-recursion}. Theorem~\ref{Thm:Q-expansions} and Corollary~\ref{Cor:Q-linearizations} provide combinatorial identities for these polynomials.

\item
Let $\nu$ be a positive measure on $\mf{R}$ with all moments finite, and $\Exp{x} = 0$, $\Exp{x^i} = m_{i-2}(\nu)$. Then $\Exp{\cdot}$ is positive on $\mf{C}_0[x]$, but in general has no positive extension to all of $\mf{C}[x]$. So we are in the context of part (b) of the preceding alternative. As a result, we may replace $\mc{A}_0$ by $\mc{A}_0 / I_{\Exp{\cdot}}$ throughout, as long as the operators $X(x^i)$ without the scalar part are used. This is the $q$-Kolmogorov case (see \cite{AnsQLevy}), which includes all the centered $q$-infinitely divisible measures (all of whose moments are finite).

\item
Let $\mc{A}_0 = \mf{C}_0[x]$, $\Exp{x^i} = t \delta_{i, 2} = t m_{i-2}(\delta_0)$. Note that this functional is positive on $\mf{C}_0[x]$, but not on $\mf{C}[x]$. Then
\[
I_{\Exp{\cdot}} = \set{P \in \mf{C}_0[x]: P'(0) = 0}
\]
is an ideal, and condition \eqref{Gaussian} is satisfied. This is the $q$-Gaussian case. In $\mc{A}_0 / I_{\Exp{\cdot}}$, $x^i = 0$ for $i \geq 2$. It follows that the $q$-Kailath-Segall polynomials in this case can also be considered as single-variable polynomials. Specifically, they are the continuous $q$-Hermite polynomials of Example~\ref{Ex:q-Hermite}. Theorem~\ref{Thm:Q-expansions} and Corollary~\ref{Cor:Q-linearizations} provide combinatorial identities for these polynomials as well.
\end{enumerate}
\end{Ex}

\subsection{Differences from the $q=0,1$ cases}
Binomial families of polynomials are Sheffer families for $t=0$. In particular, for any Appell polynomials the binomial family is always $\set{x^n}$.

\br
In \cite{AnsMeixner}, we investigated $q$-analogs of binomial families, which have generating functions
\begin{equation}
\label{Old-binomial}
\exp_q(U(z) x) = \prod_{k=0}^\infty \frac{1}{1 - q^k U(z) x}.
\end{equation}
We showed (Proposition~18) that the theory of such binomial families for various values of $q$ is exactly parallel to the theory considered by Rota et al. for $q=1$. In particular, the lowering operator (in the sense of Theorem~\ref{Thm:Characterization}) for such a sequence is $U^{-1}(D_q)$. In fact, these results were surely known before.

\br
However, the following is one difference between the classical and free Sheffer (or even Appell) families. For a classical (single variable) Sheffer family $\set{P_n(x,t)}$ with generating functions $f(z)^t e^{x U(z)}$,
\begin{equation}
\label{Binomial1}
U^{-1}(\partial_x) P_n(x,t) = P_{n-1}(x,t).
\end{equation}
However, for a corresponding free Sheffer family with generating function $\frac{1}{1 + f(z) t - U(z) x}$, the corresponding property is
\[
U^{-1}(\partial_x) P_n(x,t) = \sum_{k=0}^{n-1} P_k(x,t) P_{n-k-1}(x,t).
\]
In other words, the lowering operator for the classical family is $U^{-1}(\partial_x)$ independently of $t$. In the free case, for $t=0$ the lowering operator is $U^{-1}(D_0)$, but it changes for $t>0$.

\br
For $q \neq 0, 1$, the binomial families corresponding to the $q$-Sheffer sequences~\eqref{q-Sheffer} of this paper, with generating functions
\[
\prod_{k=0}^\infty \frac{1}{1 - U(q^k z) x},
\]
are different from the binomial families~\eqref{Old-binomial}. For binomial families of this paper, the lowering operator need not be a function of $D_q$ (which is the lowering operator for the family $\set{x^n}$). For example, the binomial family for the continuous big $q$-Hermite polynomials is
\[
\set{P_n = \prod_{k=0}^{n-1} (x - [k]_q)},
\]
and the lowering operator for this family is
\[
f \mapsto \frac{f(x) - f(1 + q x)}{x - (1 + q x)}.
\]
This operator does not commute with $D_q$.

\br
Another property which still holds in the free case is that any free Sheffer family can be expressed as a linear combination of the corresponding Appell family \cite[Lemma 1]{AnsMeixner}. This property does not hold in the $q$ case: see the appendix.

\appendix

\section{$q$-Appell polynomials are not martingale polynomials for $q$-L\'{e}vy processes}
\noindent
All the calculations below are performed with Maple~$7$.

\subsection{Relation between $q$-Appell and $q$-Sheffer polynomials}
Let $\set{P_n(x,t)}$ be the centered continuous big $q$-Hermite polynomials, and $\set{A_n(x,t)}$ the $q$-Appell polynomials for $r_1=0$, $r_k=t$ for $k>1$. Note that $\set{P_n}$ are $q$-Sheffer for this cumulant sequence. Then
\[
\begin{split}
P_4(x,t)
&- \Bigl(A_4(x,t) - (3 + 2q + q^2)A_3(x,t) + (3 + 4q + 3q^2 + q^3)A_2(x,t)  \\
&\qquad - (1 + 2q + 2q^2 + q^3)A_1(x,t) \Bigr) \\
&= t [x q (q - 1)].
\end{split}
\]
It follows that $q$-Sheffer polynomials need not be a $t$-independent linear combination of the corresponding $q$-Appell polynomials. In particular, since $\set{P_n(x,t)}$ are known to be martingale polynomials for the $q$-Poisson process \cite{AnsLinear}, it follows that the $q$-Appell polynomials, in this case, are not.

\subsection{Generic case}
Take $\mc{A}_0$ to be the polynomial algebra generated by symbols $\set{x(t): t \in \mf{R}_+}$, and the functional on it be given by
\[
\Exp{x(t_1) x(t_2) \ldots x(t_k)} = (\min_i t_i) r_k.
\]
Denote $X(t) = X(x(t))$ as in the case (b) of Section~\ref{Section:KS-Fock-space}. Then $\set{X(t)}$ is a $q$-L\'{e}vy process in the sense of \cite{AnsQCum}. For each fixed $t$, define $W_{\vec{u}}(t)$ using the recursion relation~\eqref{q-KS} involving only the $X((x(t))^k)$, for that specific $t$. The conditional expectations onto the subalgebras generated by the $q$-L\'{e}vy process $\set{X(t)}$ are determined by
\[
\Cexp{s}{W_{\vec{u}}(t)} = W_{\vec{u}}(s).
\]
In particular, for each $\vec{u}$ the Kailath-Segall family $W_{\vec{u}}(t)$ is a martingale. In a number of situations, these processes also have single-variable martingale polynomials families. That is, for each $n$ there is a family of polynomials $\set{P(x,t)}_{t \in [0, \infty)}$ of degree $n$ in $x$ such that
\begin{equation}
\label{Cond-exp}
\Cexp{s}{P(X(t), t)} = P(X(s), s).
\end{equation}
This is the case for $q=1$, $q=0$, if the process is a $q$-Brownian motion, and if the process is a $q$-Poisson process. In all of these cases, from the existence of such martingale polynomials one can deduce the Markov property for the corresponding process \cite{BKSQGauss,AnsLinear}.

\br
We show that generically, there is no degree $5$ polynomial which is a martingale for such a $q$-L\'{e}vy process. This is a strong indication that general $q$-L\'{e}vy processes do not have the Markov property. For a fixed $t$, explicit Maple calculations (see the author's web page) show that there is a monic degree $5$ polynomial $P_5(x)$, whose coefficients depend only on $t$, $q$, and the free cumulants of the process, such that
\[
\begin{split}
& P_5(X(t)) + q^2 (1-q) t \bigl(r_3 W_2(t) - r_2 W_3(t)\bigr) \\
&\quad = \text{time-independent linear combination of the Kailath-Segall polynomials}.
\end{split}
\]
Since a linear combination of Kailath-Segall polynomials is a martingale, the expression above is a martingale. Suppose there is a degree $5$ martingale polynomial in $x$. Further explicit calculations show that there definitely are martingale polynomials of degrees $4$ and less. It follows from equation~\eqref{Cond-exp} that the leading coefficients of such polynomials have to be independent of $t$, so we may assume them to be monic. By subtracting from $P_5$ the degree $5$ martingale polynomial, we obtain a martingale
\[
P_4(X(t)) + q^2 (1-q) t \bigl(r_3 W_2(t) - r_2 W_3(t)\bigr),
\]
where the polynomial $P_4$ has degree at most $4$. $P_4(X(t)) \cdot 1$ contains a term of degree $\deg P_4$, while $(r_3 W_2(t) - r_2 W_3(t)) \cdot 1 = (r_3 x(t)^2 - r_2 x(t)^3)$, so if $ \deg P_4 > 1$, $P_4$ also has to have a time-independent leading coefficient. Proceeding in this fashion, we obtain a martingale of the form
\[
Y(t) = a(t,q,r) X(t) + b(t,q,r) + q^2 (1-q) t (r_3 W_2(t) - r_2 W_3(t))
\]
Since
\[
Y(t) \cdot 1 = a(t,q,r) x(t) + b(t,q,r)  + q^2 (1-q) t (r_3 x(t)^2 - r_2 x(t)^3),
\]
$b$ is independent of $t$, and may be taken to be zero. First assume $a \equiv 0$. This can occur only in the following cases.
\begin{enumerate}
\item
$q=1$. This is the classical L\'{e}vy process case.
\item
$q=0$. This is the free L\'{e}vy process case.
\item
$r_3 W_2(t) = r_2 W_3(t)$. Suppose that condition~\eqref{Cstar} holds. From the ``creation'' part of the representation of the operator $X(t)$ in Section~\ref{Section:KS-Fock-space}, it follows that $r_3 x(t)^2 - r_2 x(t)^3 \in I_{\Exp{\cdot}}$, and in particular that
\[
\Exp{x(t)^k (r_3 x(t)^2 - r_2 x(t)^3)} = t (r_3 r_{k+2} - r_2 r_{k+3}) = 0
\]
for $k \geq 0$. This is the case only when $r_k = \alpha^{k-2}$ for $k \geq 2$ (we may take $r_1 = 0$, $r_2 = 1$). Then $W_k(t) = \alpha^{k-2} X(t)$. These are $q$-Poisson processes with step $\alpha$, and the degenerate case $\alpha = 0$ gives the $q$-Brownian motion.
\end{enumerate}
For general $a$, and $q \neq 0, 1$,
\[
\Cexp{s}{Y(t)} \cdot 1 = a(t,q,r) x(s) + q^2 (1-q) t (r_3 x(s)^2 - r_2 x(s)^3),
\]
and the martingale condition implies
\[
(a(t,q,r) - a(s,q,r)) x(s) + q^2 (1-q) (t-s) (r_3 x(s)^2 - r_2 x(s)^3) = 0.
\]
As above, this implies
\[
(a(t,q,r) - a(s,q,r)) r_{k+1} + q^2 (1-q) (t-s) (r_3 r_{k+2} - r_2 r_{k+3}) = 0
\]
for $k \geq 0$. By fixing $s$ and varying $t$, it follows that $a(t,q,r)$ is linear in $t$, and may be taken to be $a(q,r) t$. With $r_2 = 1$, $r_3 = \alpha$, this means that
\[
r_{k+3} = \alpha r_{k+2} + \frac{a}{q^2 (1-q)} r_{k+1}.
\]
For general $a, \alpha$, this says that $\set{r_k}$ is a sum of two geometric sequences, corresponding to the $q$-L\'{e}vy measure being supported at two points.

\begin{Cor}
The Markov processes of \cite{Bryc-Meixner} are not $q$-L\'{e}vy processes. 
\end{Cor}

\begin{proof}
Proposition~3.3 of \cite{Bryc-Meixner} shows that the processes of that paper have martingale polynomials.
\end{proof}


\providecommand{\bysame}{\leavevmode\hbox to3em{\hrulefill}\thinspace}
\providecommand{\MR}{\relax\ifhmode\unskip\space\fi MR }
\providecommand{\MRhref}[2]{%
  \href{http://www.ams.org/mathscinet-getitem?mr=#1}{#2}
}
\providecommand{\href}[2]{#2}

\end{document}